\documentclass[onefignum,onetabnum,reqno]{siamonline190516}

\usepackage{enumerate,subfig}
\usepackage{amsfonts}
\usepackage{graphicx}
\usepackage{epstopdf,wrapfig,cancel}

\usepackage{color}
\usepackage{amsopn}
\usepackage{tikz}
\usetikzlibrary{arrows}

\usepackage{verbatim}
\ifpdf
  \DeclareGraphicsExtensions{.eps,.pdf,.png,.jpg}
\else
  \DeclareGraphicsExtensions{.eps}
\fi

\usepackage{enumitem}
\setlist[enumerate]{leftmargin=.5in}
\setlist[itemize]{leftmargin=.5in}

\newcommand{\nc}{\newcommand}
\nc{\N}{\mathbb{N}}
\nc{\Z}{\mathbb{Z}}
\nc{\D}{\mathbb{D}}
\nc{\Q}{\mathbb{Q}}
\nc{\R}{\mathbb{R}}
\nc{\C}{\mathbb{C}}
\nc{\T}{\mathbb{T}}
\nc{\cX}{\mathcal{X}}
\nc{\cY}{\mathcal{Y}}
\nc{\cL}{\mathcal{L}}
\nc{\cH}{\mathcal{H}}
\nc{\cT}{\mathcal{T}}
\nc{\cU}{\mathcal{U}}
\nc{\cV}{\mathcal{V}}
\nc{\taub}{\boldsymbol{\tau}}
\nc{\btheta}{\boldsymbol{\theta}}
\nc{\Sph}{\mathbb{S}^2}
\nc{\tld}[1]{\tilde{#1}}
\nc{\wtld}[1]{\widetilde{#1}}
\nc{\hu}{\hat{u}}
\nc{\wh}[1]{\widehat{#1}}
\nc{\Pb}{\mathbb{P}}
\nc{\Sb}{\mathbb{S}}
\nc{\Fbf}{\textbf{F}}
\nc{\Gbf}{\textbf{G}}
\nc{\Lbf}{\textbf{L}}
\nc{\Nbf}{\textbf{N}}
\nc{\Ibf}{\textbf{I}}
\nc{\Dbf}{\textbf{D}} 
\nc{\Tbf}{\textbf{T}}
\nc{\mbf}{\textbf{m}} 
\nc{\Rbf}{\textbf{R}}   
\nc{\xb}{\textbf{x}}
\nc{\sbb}{\textbf{s}}
\nc{\ub}{\textbf{u}}
\nc{\tb}{\textbf{t}}
\nc{\rb}{\textbf{r}}
\nc{\gb}{\textbf{g}}
\nc{\hb}{\textbf{h}}
\nc{\vb}{\textbf{v}}
\nc{\eb}{\textbf{e}}
\nc{\wb}{\textbf{w}}
\nc{\bb}{\textbf{b}}
\nc{\alphab}{\boldsymbol{\alpha}}
\nc{\betab}{\boldsymbol{\beta}}
\nc{\etab}{\boldsymbol{\eta}}
\nc{\epsb}{\boldsymbol{\varepsilon}}
\nc{\phib}{\boldsymbol{\varphi}}
\nc{\fb}{\boldsymbol{f}}
\nc{\xib}{\boldsymbol{\xi}}
\nc{\gammab}{\boldsymbol{\gamma}}
\nc{\zetab}{\boldsymbol{\zeta}}
\nc{\nub}{\boldsymbol{\nu}}
\nc{\thetab}{\boldsymbol{\theta}}
\nc{\cb}{\boldsymbol{c}}
\nc{\Wb}{\boldsymbol{W}}
\nc{\Ab}{\boldsymbol{A}}
\nc{\Ib}{\boldsymbol{I}}
\nc{\Xb}{\boldsymbol{X}}
\nc{\Kb}{\boldsymbol{K}}
\nc{\Nb}{\boldsymbol{N}}
\nc{\yb}{\textbf{y}}
\nc{\zb}{\textbf{z}}
\nc{\kb}{\textbf{k}}
\nc{\qb}{\textbf{q}}
\nc{\ph}{\varphi}

\DeclareMathOperator{\sign}{sign}

\newsiamremark{remark}{Remark}
\newsiamthm{newthm}{Theorem}
\newsiamremark{newdef}{Definition}
\newsiamthm{cor}{Corollary}
\newsiamthm{prop}{Proposition}
\newsiamremark{hypothesis}{Hypothesis}
\crefname{hypothesis}{Hypothesis}{Hypotheses}
\newsiamthm{claim}{Claim}

\newcommand{\abs}[1]{\left\lvert#1\right\rvert}
\newcommand{\norm}[1]{\left\lVert#1\right\rVert}

\headers{Discovery of Dynamics}{Keller and Du.}

\title{Discovery of Dynamics using Linear Multistep Methods\thanks{Submitted to \emph{SIAM Journal on Numerical Analysis}.
\funding{This research was supported by the NSF Graduate Research Fellowship Grant DGE-1644869 for R.K., and by
NSF DMS-1719699, DMS-2012562, CCF-1704833  and ARO MURI Grant W911NF-15-1-0562
for Q.D.}}}

\author{Rachael Keller\thanks{Department of Applied Physics and Applied Mathematics, Columbia University, New York, NY
  R.K. \email{rachael.keller@columbia.edu}, Q.D. \email{qd2125@columbia.edu})}
\and Qiang Du\footnotemark[1]
}

\ifpdf
\hypersetup{
  pdftitle={Discovery of Dynamics using Linear Multistep Methods},
  pdfauthor={Keller, R, and Q. Du}
}
\fi

\begin{document}

\maketitle

\begin{abstract}
Linear multistep methods (LMMs) are popular time discretization techniques for the numerical solution of differential equations.  Traditionally they are applied to solve  for the state given the dynamics (the forward problem), but here we consider their application for learning the dynamics given the state (the inverse problem). This repurposing of LMMs is largely motivated by growing interest in data-driven modeling of dynamics, but the behavior and analysis of LMMs for discovery turn out to be significantly different from the well-known, existing theory for the forward problem. Assuming {a highly idealized setting of being given the exact state with a zero residual of the discrete dynamics,} we establish for the first time a rigorous framework based on refined notions of consistency and stability to yield convergence using LMMs for discovery. 
 When applying these concepts to three popular $M-$step LMMs, the Adams-Bashforth, Adams-Moulton, and Backwards Differentiation Formula  schemes, 
  the new theory suggests that Adams-Bashforth {for  $M$ ranging from $1$ and  $6$, Adams-Moulton for $M=0$ and $M=1$, and  Backwards Differentiation Formula for all positive $M$}  are convergent, and, otherwise, the methods are not convergent in general. In addition, we provide numerical experiments to both motivate and substantiate our theoretical analysis.
\end{abstract}



\section{Introduction}
In this work, we focus on {developing a new numerical analysis framework for}  the \emph{discovery} of dynamical systems with given states, where finitely many discrete measurements are used to approximately recover the unknown dynamical system --  a \emph{data-driven} discovery of dynamics \cite{brunton2019data,lipson09}. {Data-driven discovery of dynamical systems is experiencing a renaissance as costs of sensors, data storage, and computational resources has decreased \cite{brunton17}. {Meanwhile}, 
 advancements in the fields of machine learning and data science \cite{goodfellow2016deep,jordan15,krizhevsky12,lecun15,shalev2014} have brought in renewed vigor and enabled expansive view to this field.  At the same time, the growth of data-driven
discovery of dynamical systems has also led to a new solution method and model reduction approach to study  multiscale and high dimensional complex problems. For more discussions, we refer to  works such as 
\cite{bhattacharya20,brunton16,gulgec19,han2018solving,kang19,kevrekidis2016kernel,khoo17,long19,lu19,ma2018model,pan2018data,qian2020lift,qin19,raissi18dl,mnn,raissi18gp, pinn, rudy19, sun19,perdikaris19,wang2016deep,williams2015data,wu19,webster18,zhu18}.

\subsection{Motivation:  Data-driven discovery of dynamical systems via  linear multistep methods}
{In this work, we consider using linear multistep methods (LMMs) to discover unspecified dynamics given the state at equidistant time steps and contribute to the fundamental theory of using LMMs for data-driven discovery.  Historically, LMMs have been developed as popular schemes for numerically integrating known dynamic systems \cite{goldstine12},  with well-established mathematical theory in the last century \cite{atkinson11,dahlquist63,gautschi,henrici62,suli03}.}  
Recent works {combine} the classical numerical technique of linear multistep methods with neural networks for dynamics discovery \cite{mnn, perdikaris19,webster18}.

\begin{figure}[H]\centering
\subfloat[AM with 256 Width Network]{\includegraphics[width=0.45\columnwidth]{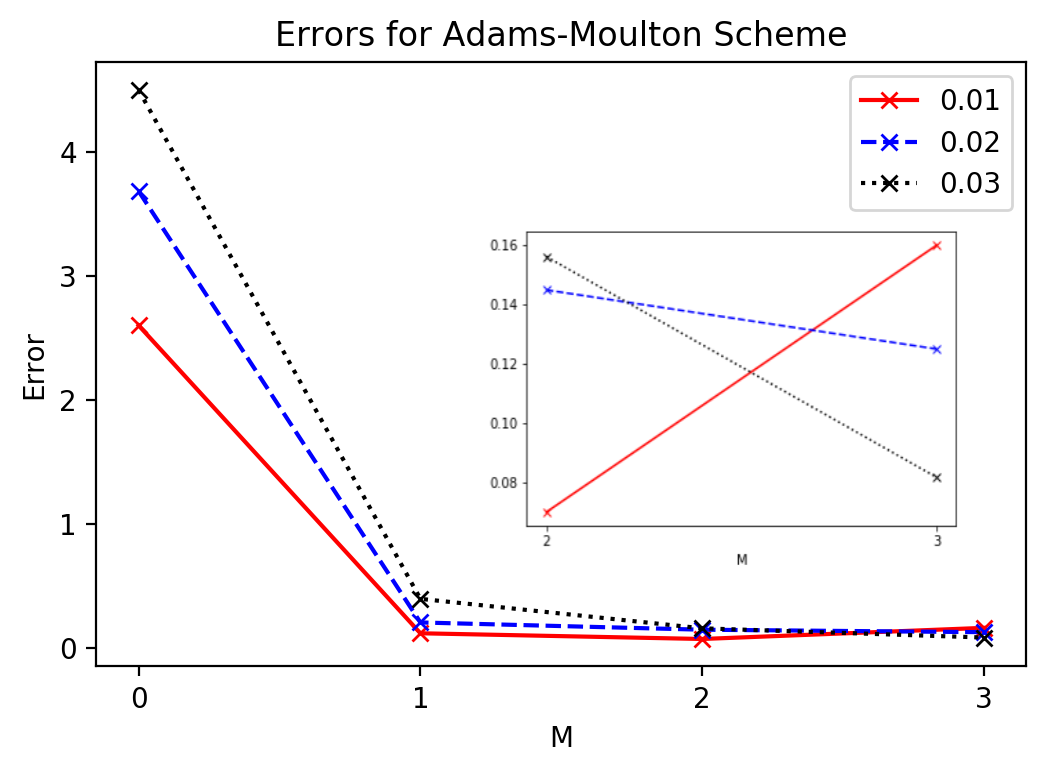}}
\subfloat[AM with 512 Width Network]{\includegraphics[width=0.45\columnwidth]{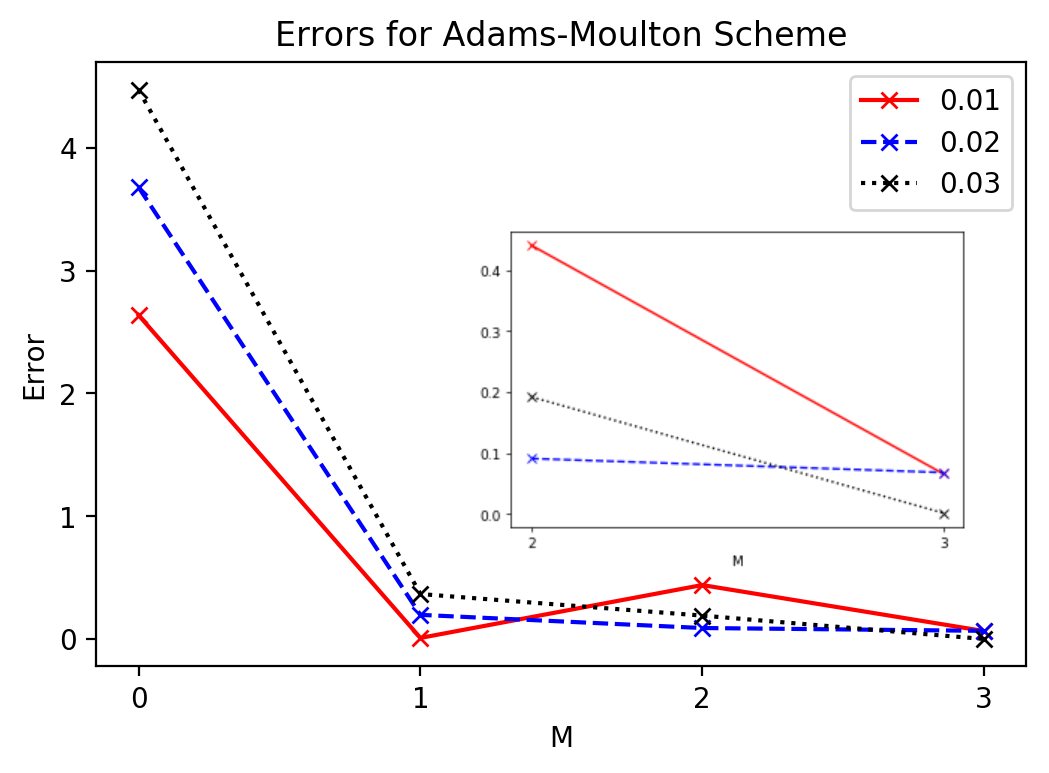}\label{fig:am512}}\\
\subfloat[AB with 256 Width Network]{\includegraphics[width=0.45\columnwidth]{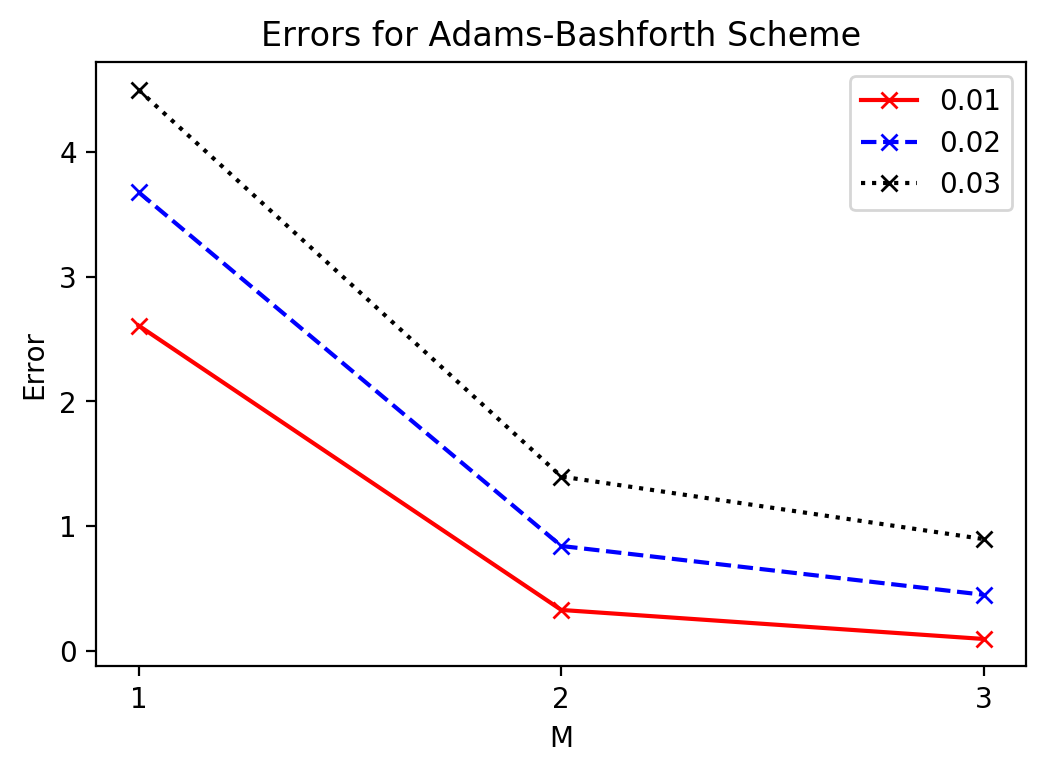}}
\subfloat[BDF with 256 Width Network]{\includegraphics[width=0.45\columnwidth]{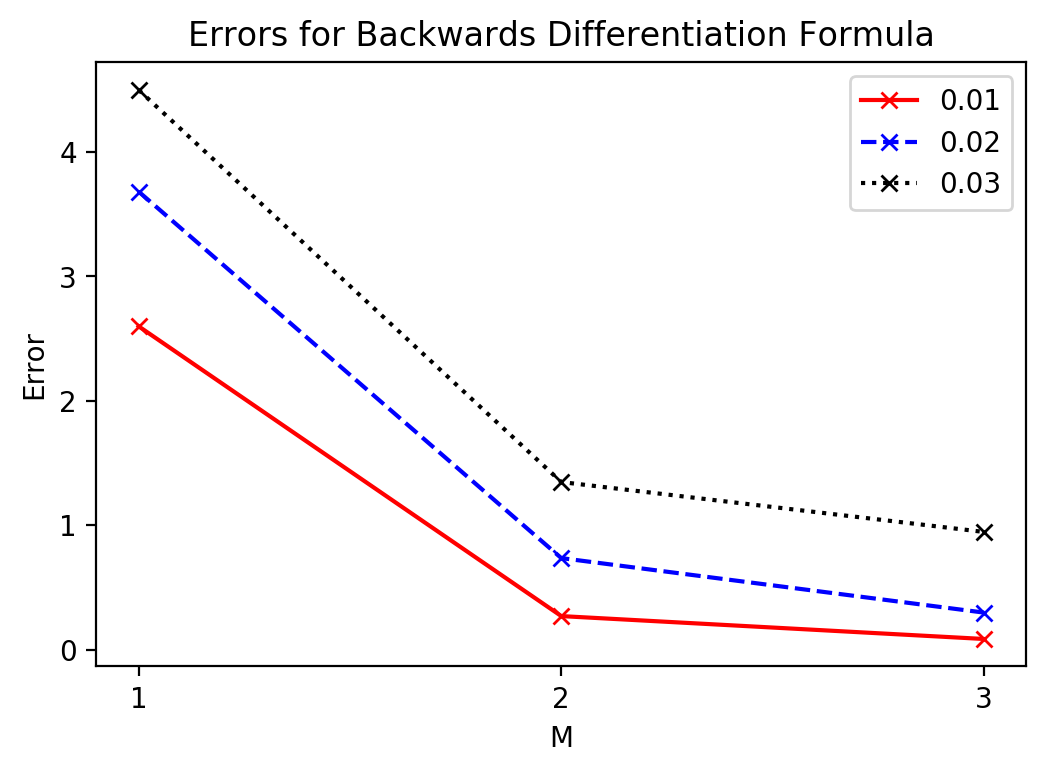}}
\caption{Absolute  $\ell_2$-errors for the first coordinate of the 2D Damped Cubic System \eqref{model1} on $t \in [0,5]$ with varying time mesh {size} $h = 0.01, 0.02, 0.03$, using a single hidden layer neural network with $\tanh$ activation function, as used in \cite{mnn}, after a fixed number of training iterations for each $M$.}\label{fig:ab-am-bdf0}
\end{figure}

Coined ``LMNet,'' LMMs are combined with neural networks for discovery of dynamics in \cite{mnn,perdikaris19,webster18}.  Figure \ref{fig:ab-am-bdf0} shows the absolute errors associated with learning $\fb$ for a nonlinearly-damped, 2D cubic oscillator \eqref{model1} using neural networks with three representative schemes of LMMs -- Adams-Moulton (AM), Adams-Bashforth (AB), and Backwards Differentiation Formula (BDF). These results are generated using the code repository built for \cite{mnn}; reported are the errors of the dynamics rather than the integrated dynamics, which are shown in \cite{mnn}.  For solving differential equations with smooth solutions, increasing $M$ corresponds to higher accuracy if the scheme is also stable. The AM scheme is an example of such a method; hence, the perplexing behavior in the errors of AM as observed in { \cite{mnn,webster18} (see Tables 1 and 2 of \cite{mnn} and Table 1 of \cite{webster18}). {As $M$ increases and $h$ decreases, the errors do not decrease. Further, as we expand the width, thereby increasing the expressibility of the network, the scheme still does not exhibit stable behavior. On the other hand,} as shown in Figure \ref{fig:ab-am-bdf0}, the AB and BDF methods with a fixed network size of 256 show a trend of convergence as $M$ and the mesh size $h$ decrease, while the AM methods show erratic behavior for the same width, persistent even with more expressibility of the network by widening the hidden layer (Figure \ref{fig:am512}). Since AM is a stable method as a time integrator, these findings warrant further investigation.} {Indeed, it has also been observed by others that increased resolution does not necessarily imply better neural network representation and prediction without a mathematically sound formulation of the learning problem \cite{bhattacharya20}. While there are many contributing factors such as the neural network structure and size as well as the training process,} 
it is the goal of this paper to investigate these findings and provide a theoretical explanation of the phenomena. 

{To begin, we pose the problem of discovery of dynamics.} 
 {In contrast to numerically integrating dynamics to learn the state, as many classical numerical methods do, this study focuses on learning the dynamics given the state. Dynamics discovery may be viewed as an inverse problem to the forward problem of classical numerical integration. is using a classical numerical technique for the inverse problem. Well-studied for the forward problem, LMMs in this inverse setting raises questions of classical notions of consistency, stability, and convergence. We seek in this work to investigate if the classical theory  for LMMs  as time integrators to solve the forward problem has an analog or counterpart in solving the inverse problem of learning dynamics.}
To initiate studies in this direction, we introduce a systematic framework for the numerical analysis of discovery of dynamics using LMMs.  Our new framework is rooted in the classical theory for LMMs as numerical integrators of differential equations, but it adopts new stability and convergence criteria due to the inverse  nature of using discrete time integrators for dynamics discovery.}  Consequently, it draws different conclusions regarding convergence in stark contrast to the conventional wisdom. The stability properties of particular schemes depart from the { traditional numerical differential equation viewpoint, and some methods that are stable for the forward problem do not retain the property for the inverse problem dynamics discovery.} Our theory is able to explain the unusual phenomena as reported in Figure \ref{fig:ab-am-bdf0}  and 
lays a rigorous foundation for further elucidating the effect of neural networks on dynamics discovery via LMMs {through follow-up studies}. Therefore, this helps the scientific community broadly in our goal of making machine learning more transparent, explainable, stable and trustworthy.

\subsection{Summary of Results}

We present a framework in Section \ref{sec:discovery-intro} consisting of nuanced notions of consistency and stability to handle unique challenges presented by using LMMs for discovery. These concepts are then combined to prove convergence. {A set of algebraic criteria is developed to check for the consistency and stability, and thus convergence, of LMMs for dynamics discovery.} With this foundation, in Theorems  \ref{thm:consist-dynamics} and \ref{thm:scheme-stability} {, } we outline consistency and stability properties of the Adams-Bashforth, Adams-Moulton, and Backwards Differentiation Formula schemes, and consequentially, Corollary \ref{cor:conv}, their convergence guarantees.

\subsection{Outline}

This paper is organized as follows. In Section \ref{sec:prob-intro} we briefly review LMMs and their theory for solving ordinary differential equations, including the standard notions in numerical analysis of truncation error,  consistency, stability, and convergence, along with an algebraic root condition for stability. In Section \ref{sec:discovery-intro} we frame the problem of discovery using LMMs and develop nuanced versions of consistency and stability for discovery. 
 In particular, in Section \ref{sec:trunc-error}, we discuss how truncation error for discovery is inherited from the forward problem and  introduce a stronger notion of consistency; in Section \ref{sec:stability} we refine the traditional definition of stability 
 and the algebraic root condition, and we show {equivalent theorems connecting the root conditions and the refined notions of stability.} In Section \ref{sec:app}, the discovery framework of Section \ref{sec:discovery-intro} is applied to characterize convergence properties of the Adams-Bashforth,  Adams-Moulton, and Backwards Differentiation Formula schemes. Some discussions on the long time dynamics discovery are made in 
Section \ref{sec:long-time}.
 Then, in Section \ref{sec:numerics}, we show results of numerical experiments.  
 Finally, in Section \ref{sec:conc}, we summarize the results and discuss future directions.

\section{LMMs: Quick Review}\label{sec:prob-intro} 
In this section, we introduce notation used throughout this work and briefly review the theory of LMMs as time integrators. While LMMs are well-documented in standard textbooks for solving ordinary differential equations (see \cite{gautschi, suli03, atkinson11, henrici62}), we include the salient points to facilitate direct comparison with the new theory for the discovery of unknown dynamics developed in the next section.

\subsection{LMMs: Notation and Concepts}\label{sec:notation}
Consider the ordinary differential equation {(ODE)}
\begin{equation}\label{model0}
\frac{d}{dt}\xb(t) = \fb(\xb(t)), \ a \leq t \leq b,   \ \xb(t_0) = \xb_0,
\end{equation}
where $\xb \in C^\infty(0,\infty)^d$ and $f$ is assumed to be a Lipschitz continuous, smooth, and bounded function. Discretizing the model problem \eqref{model0}, we assume a grid on the interval $[ \, a,b\, ]$ defined to be a set of points: $a=t_0 < t_1 < \cdots < t_{N} = b$ with equidistant mesh  {$t_{n+1} - t_n=h=(b-a)/(N+1)$, $n \in\{ 0, 1, \ldots, N\}$}. Let $[ \, a,b\, ]_h$ denote this ordered set. We denote the set of grid functions $\Gamma_h[ \, a,b\, ] = \left\{ \zb \, \vert \, \zb \in \R^{{(N+1)} \times d}, \ \zb_n = \zb(t_n)  \in \R^d, t_n \in [ \, a,b\, ]_h\right\}$ \cite{gautschi}. 

{An $M$-step LMM approximates the $n^{th}$ value $\xb_n = \xb(t_n)$ in terms of the previous $M$ ($M\geq 1)$ time steps $\xb_{n-1}, \xb_{n-2}, \ldots, \xb_{n-M}$} \cite{gautschi, suli03, atkinson11, henrici62}. An $M-$step linear multistep method is given by
\begin{equation}\label{lmm000}
\begin{split}
\sum_{m=0}^M \alpha_m \xb_{n-m} &= h \sum_{m=0}^M \beta_m \fb(\xb_{n-m}), \ n = M, M+ 1, \ldots, N,
\end{split}
\end{equation}
where $\xb \in \Gamma_h[ \, a,b\, ]$, the coefficients $\alpha_m, \beta_m \in \mathbb{R}$ for $m= 0, 1, \ldots, M,$ and  $\alpha_0 \not=0.$ The function $\fb$ is assumed to be given and Lipschitz, and the LMM scheme \eqref{lmm000} defines an iterative procedure stepping forward in the independent variable $t \in [\, a,b \, ]$ to solve for $\xb(t)$ at the gridpoints.  {Associated with an $M-$step LMM are its first and second characteristic polynomials, given, respectively,  by
\begin{equation}
\label{eqn:char}
\rho(z) = \sum_{m=0}^M \alpha_{M-m} z^{m}, \quad\text{and}\quad
\sigma(z) = \sum_{m=0}^M \beta_{M-m} z^{m},
\end{equation}
 where it is assumed that $\alpha_0 \not=0$ \cite{suli03}.}

For the numerical integration of differential equations, the method \eqref{lmm000} is called explicit if $\beta_0 = 0$ and implicit otherwise \cite{gautschi, suli03, atkinson11}. Implicit methods require a nonlinear solver to the generated system of equations, whereas explicit methods do not. Existence and uniqueness of solutions in the case of implicit schemes is shown in \cite{gautschi, henrici62}. For both implicit and explicit methods, a kickstarting method for initial $M$ values must be chosen, and as such a critical component of analyzing any multistep method scheme is to understand how much errors in initial values pollute the subsequent calculations \cite{gautschi}. 
This aspect of numerical methods is called numerical stability  \cite{atkinson11}. 

Finally, {for any index set $\mathcal{S}$ with cardinality $\bar{S}$, 
we let $\norm{\zb}_1 =\sum_{i\in \mathcal{S}}  \abs{\zb_i}$ and
 $\norm{\zb}_\infty = \max_{i \in \mathcal{S}} \abs{\zb_i}$ denote the standard discrete norms for
 any vector $\zb$ naturally embedded in $\mathbb{R}^{\bar{S}\times d}$
 where $\abs{\zb_i}$ can be any vector norm of $\zb_i\in\mathbb{R}^d$.
 The same notations are used also for discrete grid functions given  either in $\Gamma_h[a,b]$ or its subsets.}

\begin{remark}
To fix ideas, we use the hat notation $\hat{\; }$ to mark grid functions generated by the discretization \eqref{lmm000}.
 In the forward problem, the state $\xb(t)$ is iteratively produced by LMMs, and hence we study $\hat{\xb}$, whereas for dynamics discovery, we study $\hat{\fb}$, see Section \ref{sec:discovery-intro}. 
\end{remark}

\subsection{The Adams Family and BDF}\label{sec:methods}
Adams-Bashforth (AB), Adams-Moulton (AM), and the Backwards Differentiation Formula (BDF) are three popular multistep method schemes that arise from a Lagrange interpolating polynomial of the state or dynamics at time $t_n$ using data from previous time steps. Without loss of generality, we consider the scalar model problem in this section; for higher dimensions, the theory need only be applied in each dimension. 
{Let $\Lambda_0 = \{ -M, -M+1, \ldots, -1, 0\}$
 and  $\Lambda_1 = \{-M, -M+1, \ldots, -1\}$.  The Lagrange interpolating polynomial  of a function $u : \R \rightarrow \R$ over the set  $\{t_{n+i}, i\in \tilde{\Lambda}\}$ 
 is the polynomial of degree $M$ for $\tilde{\Lambda}=\Lambda_0$ and degree $M-1$ 
  for $\tilde{\Lambda}=\Lambda_1$ obtained from the linear combination of basis functions}
\begin{equation}\label{eqn:lag-basis-ele}
{ \ell_{k,n}(t;  \tilde{\Lambda}) = \prod_{\substack{i\in \tilde{\Lambda} \setminus \{k\}}}\frac{t-t_{n+i}}{t_{n+k}-t_{n+i}}, \ k \in \tilde{\Lambda}, }
\end{equation}
{with  $u(t_{n+k})$ for each $k\in\tilde{\Lambda}$ as the coefficient of the linear combination.
{The M-step} Adams-Moulton (or AM-$M$)  and  Adams-Bashforth (or AB-$M$) are 
 $M$-step LMMs that arise from  interpolating  the dynamics {$f(x(t_n)),$} with
Lagrange interpolating polynomials corresponding to $\tilde{\Lambda}=\Lambda_{0}$ and
$\tilde{\Lambda}=\Lambda_1$, respectively,
and then applying the fundamental theorem of calculus on the model problem \eqref{model0}. Letting $f(t_n)$ denote $f(x(t_n))$ for brevity, we have }
\begin{align}\label{eqn:dadam00}
x(t_{n})  &\approx   x(t_{n-1}) + \int_{t_{n-1}}^{t_n} \sum_{k\in \tilde{\Lambda}} f(t_{n+k}) \ell_{{k,n}}(t; \tilde{\Lambda}){dt}\,.
\end{align}
{BDF-$M$, on the other hand, is an $M$-step LMM for $M\geq 1$ derived from interpolating the state  $\xb \in \Gamma_h[a,b]$ in \eqref{model0} directly on the lattice $\Lambda_0$, so that} 
\begin{align*}
{
\sum_{k \in \Lambda_0} x(t_{n+k}) \frac{ d\ell_{k,n} }{dt}(t_n; \Lambda_0) \approx 
\frac{d}{dt}x(t_n) = f(t_n).
}
\end{align*}
{ By the change of variables $u = (t-t_{n-1})/h$, we have a scaled Lagrange interpolating polynomial, denoted $\ell^h_k$, given by}
\begin{equation}\label{eqn:lag-basis-ele-uniform}
\ell^h_k({u}; {\tilde{\Lambda}}) = \prod_{\substack{i\in\tilde{\Lambda} \setminus \{k\}}}\frac{u-1-i}{k-i}, \ k \in \tilde{\Lambda}.
\end{equation}
With \eqref{eqn:lag-basis-ele-uniform}, the integrand  of \eqref{eqn:dadam00} may be written {independent of the time step, so that}
{
\begin{align}\label{eqn:dadam00-uniform}
x(t_{n})  &\approx   x(t_{n-1}) + \int_{0}^{1} \sum_{k\in \tilde{\Lambda}} f(t_{n+k}) \ell^h_k(u; \tilde{\Lambda}){du}.
\end{align} 
The simplified coefficients for the BDF method with uniform mesh can be obtained similarly.
 }

\subsection{Truncation Error and Consistency} {In this section, we introduce the residual and notions related to analytical error for LMMs.
The residual operator is given by \cite{gautschi}:}
\begin{equation}\label{eqn:res-forward}
{(R_h \hat{\xb})_n := \frac{1}{h} \sum_{m=0}^M \alpha_m \hat{\xb}_{n-m} - \sum_{m=0}^M \beta_m \fb(\hat{\xb}_{n-m}), \quad 
 n = M, M+1, \ldots, N,}
\end{equation}
{defined for $\hat{\xb} \in \Gamma_h[ \, a,b\, ]$.}  How accurately the discretization \eqref{lmm000} approximates the solution of \eqref{model0} is measured by the truncation error, defined below.

\begin{newdef}[Local Truncation Error \cite{suli03, atkinson11, henrici62, gautschi}]\label{def:trunc}
{Let $\xb \in \Gamma_h[a,b]$ be the exact solution of the dynamic system \eqref{model0} defined at the grid coordinates. The local truncation error $\taub_h=
((\taub_h)_M, (\taub_h)_{M+1}, \ldots, (\taub_h)_{N}) \in \R^{(N-M+1)\times d}$
 is given by}
\begin{equation}
(\taub_h)_n =(R_h {\xb})_n ,\quad \mbox{ for }\; n = M, M+1, \ldots, N. \label{eq:taubhdef}
\end{equation}
For smooth functions $\fb$ and $\xb$, we have
$$
(\taub_h)_n =\sum_{m=0}^\infty C_m  h^{m-1}\nabla_{{t}}^{m} \xb(t_n),\quad \mbox{ for }\; n = M, M+1, \ldots, N,
$$
where
\begin{align}
C_0 =\ \sum_{k=0}^M \alpha_k,   \ 
C_m = (-1)^m\left[\frac{1}{m!} \sum_{k=1}^M k^m \alpha_k + \frac{1}{(m-1)!} \sum_{k=0}^M k^{m-1}\beta_k\right], \ m = 1, 2, \ldots.\label{eq:c0m}
\end{align}
\end{newdef}
\noindent Now, we proceed to define order of error and the notion of consistency.

\begin{newdef}[Order of Error \cite{gautschi}]\label{def:err} A linear multistep method has error order of $p$ if $\norm{\taub_h}_\infty = \mathcal{O}(h^p)$ as $h \rightarrow 0$ and admits a \emph{principal error function} $\eb(t) \in C[ \, a,b\, ]$ provided
\[
\eb(t) \not=\boldsymbol{0} \text{ and } (\taub_h)_n = \eb(t_n) h^p + \mathcal{O}(h^{p+1}) \text{ as } h \rightarrow 0,
\]
or simply, $\norm{\taub_h - h^p \eb}_\infty = \mathcal{O}(h^{p+1}).$
\end{newdef}

\begin{newdef}[Consistency \cite{gautschi}]\label{def:consistency} A linear multistep method is consistent with the differential equation provided $\norm{\taub_h}_\infty \rightarrow 0$ as $h \rightarrow 0$.
\end{newdef}

 The Adams family and BDF are consistent in the sense of Definition \ref{def:consistency}. Moreover, the local truncation error associated with the $M-$step AB and BDF schemes are $\mathcal{O}(h^M),$ whereas for the $M-$step AM, the local truncation error is $\mathcal{O}(h^{M+1})$  \cite{suli03, atkinson11}. 

{It is well-known that consistency can be formulated algebraically in terms of the
characteristic polynomials \cite{dahlquist56}. In particular, the consistency condition, i.e., $C_0=C_1=0$
in \eqref{eq:c0m}, is equivalent to $\rho(1)=0$ and
$\rho'(1)=\sigma(1)$. Moreover,  the truncation error is order $k$ if 
\begin{align} \label{eq:order}
\rho(e^z)-z \sigma(e^z)=O(z^{k+1}),\quad\text{ as }\quad z\to 0.
\end{align}}
\subsection{Stability and the Root Condition}\label{sec:stab}
In this section, we review definitions of stability and the root condition for LMMs.  Stability is defined as follows.
\begin{newdef}[Stability \cite{gautschi}]\label{def:stable}
A linear $M-$step method for the ordinary differential equation $\dot{\xb} = \fb(t,\xb(t))$ is called stable on $[ \, a,b\, ]$ provided there exists a constant $K$ not depending on $h$ such that, for any two grid functions $\ub, \vb \in \Gamma_h[ \, a,b\, ]$, we have for all $h$ sufficiently small
\[
\norm{\ub - \vb}_\infty \leq K \left(\max_{0 \leq i \leq M-1}\norm{\boldsymbol{u}_i - \boldsymbol{v}_i} + \norm{R_h \boldsymbol{u} - R_h \boldsymbol{v}}_\infty  \right).
\]
\end{newdef}

The characteristic polynomials defined in \eqref{eqn:char}
may be used to determine the stability of a linear multistep method via the root condition.

 \begin{newdef}[Algebraic Root Condition \cite{suli03, gautschi} ]\label{def:root-condition}
 A polynomial satisfies the root condition provided the roots of the polynomial do not exceed magnitude 1, and those of magnitude 1 are simple. 
\end{newdef}

The following theorem states the equivalence between the stability and the root condition.

\begin{newthm}[Stability and the Root Condition \cite{suli03, gautschi}] 
\label{thm:root-state}
A linear multistep method is stable if and only if its first characteristic polynomial  $\rho(z)$ satisfies the algebraic root condition given by Definition \ref{def:root-condition}.
\end{newthm}

{Note that all AB and AM schemes satisfy the root condition and are stable by Definition \ref{def:stable}, whereas  BDF-$M$ satisfies the root condition and is stable only for $1\leq M\leq 6$ \cite{henrici62}.}

\subsection{Convergence}
{Finally, we introduce the definition of convergence for LMMs and the celebrated equivalence theorem for determining it.}

\begin{newdef}[Convergence \cite{gautschi}]\label{def:convergence}
Consider the initial value problem \eqref{model0} and a fixed linear multistep method defined by \eqref{lmm000}. Let {$\hat{\xb}  \in \Gamma_h[ \, a,b\, ]$} be the grid function obtained by applying 
\eqref{lmm000} on a uniform, real-valued grid of $[ \, a,b\, ]$ with mesh size $h$, and let {$\xb \in \Gamma_h[ \, a,b\, ]$} be the exact solution of \eqref{model0} at the grid points. The linear multistep method is said to converge on $[ \, a,b\, ]$ if 
\[ \norm{\xb-\hat{\xb}}_\infty \rightarrow 0 \;\mbox{ as }\; h\rightarrow 0\;\mbox{ whenever }\; \max_{0\leq k \leq M-1} \norm{{\hat{\xb}_k - \xb(t_k)}}_\infty \rightarrow 0.
\]
\end{newdef}

\noindent With Definition \ref{def:convergence}, one can obtain the Dahlquist Equivalence Theorem, Theorem \ref{thm:equiv} \cite{suli03}.

\begin{newthm}[Equivalence Theorem \cite{gautschi}]\label{thm:equiv} The multistep method \eqref{lmm000} converges in the sense of Definition \ref{def:convergence} for all Lipschitz $\fb$ if and only if it is consistent and stable.
\end{newthm}

{From the Equivalence Theorem, it can be shown that the order of the error $\norm{\xb-\hat\xb}_\infty$ is the same order as the truncation error (Definition \ref{def:err}) and thus the order of approximation, provided the initial error $\max_{0\leq k \leq M-1}\abs{\hat\xb_k - \xb(t_k)}$ is also of the same order.} 

In this work, we develop an analogous theory for multistep methods modifying these theorems to deal with the discovery of dynamics rather than solving the differential equation. In particular, we show how the second characteristic polynomial is determinant of stability for discovery and whether the Adams family and BDF are stable or not.

\section{Discovery of Dynamics}\label{sec:discovery-intro}
In this study, we consider a \emph{data-driven} technique to solve for the dynamics $\fb$ given information on the state $\xb$ at equidistant time steps \cite{mnn}. First, we introduce the problem and then discuss notions of consistency, stability, and convergence. We now proceed to define the problem of LMMs for discovery.%

\subsection{Problem Definition}
Following earlier discussions,   we are concerned with the initial value problem  \eqref{model0}.
{In this section and the next,} multivariate functions representing the continuum models
are denoted by scalar notations, i.e.,  $f=f(x)$ and $x=x(t)$, so that boldface symbols can be
reserved for vectors corresponding to discrete forms of dynamics, which should be clear in context without ambiguity.
{The  task  of learning is to produce a function to approximately represent the dynamics, $f=f(x)$, based on a set of
observed states, that conforms with the discrete dynamics described by a linear multistep method.
 In practice,  one often encounters situations with only partial (incomplete) data or data containing observation errors and uncertainties; these complications are typical for inverse problems. When combined with deep networks, 
the approximation is produced by a network in a learned parametrized form, which introduces further approximations as well as implicit regularizations. }

{
As the first step to develop a rigorous numerical analysis framework, we consider a very idealized setting in this work by assuming
that (A1) a complete set of exact values of the state, {$\{\xb_n=x(t_n)\}_{n=0}^N$,}  given at equally distributed{, ordered grid points $\{t_n\}_{n=0}^N$;} (A2) the neural networks (or the underlying function classes used to represent the dynamics) have sufficient approximation capability to produce zero residual for the discrete dynamical system; (A3) approximated values of the exact dynamics for some observed initial states are available.}

{Although the assumptions make the situation very idealized, the study is a very constructive step towards the understanding of the mathematical and computational issues related to the data-driven modeling using neural networks and discretized forms of the unknown dynamics, which are the focuses of our ongoing work. The findings made here shed light on future studies of similar issues under more realistic conditions, as discussed in Section~\ref{sec:connect} and further in Section~\ref{sec:conc}.
}
{Under the assumptions {(A1), (A2), and (A3)} stated above, the procedure of learning dynamics can be described as follows. Given $\xb_n={x}(t_n)$ for $0\leq n\leq N$ and $\hat{\fb}_{i}$ as suitable approximations of ${f}(\xb_i)$ for $i$ 
in a suitable subset of  $\{0 \leq i \leq M-1$\},  we have zero residuals for the discrete dynamics based on the LMM discretization for $t_n$ with $n=M,\ldots, N$, i.e,}
\[ 
\sum_{m=0}^M \beta_m \hat{\fb}_{n-m} = \frac{1}{h}\sum_{m=0}^M \alpha_m \xb_{n-m} ,\quad n=M , M+1, \ldots, N.
\]
Indeed, {the above system for $\hat{\fb}$} is simply  \eqref{lmm000} rewritten for learning the dynamics rather than the state. 
{To help with later discussions, we let $N_M=N-M+1$ denote the number of linear equations in the system. Given that the values of $\{\beta_m\}_{m=0}^M$ affect
the structure of the resulting system, we
 let $m_0$ and $M_0$ be the smallest and the largest index, respectively, among those $m$'s satisfying $\beta_m\neq 0$, i.e.
$$\beta_m= 0\; \text{ for any $m$ with either } m<m_0 \text{ or } m> M_0, \, \text{ while } \beta_{m_0}\neq 0 \text{ and }
\beta_{M_0}\neq 0.$$ We collect the ordered coefficients of the LMM scheme in the vectors $\alphab = (\alpha_0, \alpha_1, \ldots, \alpha_M)$ and $\betab = (\beta_{m_0}, \beta_{m_0+1}, \ldots, \beta_{M_0}).$ 
 The  system for $\hat{\fb}$ {in this reduced notation} is then
\begin{equation}
\sum_{m=m_0}^{M_0}  \beta_m \hat{\fb}_{n-m}= \frac{1}{h}\sum_{m=0}^M \alpha_m \xb_{n-m} ,\quad n=M , M+1, \ldots, N.
 \label{dyna-eqn}
\end{equation}
}

{For brevity, we introduce the index sets $\mathcal{I} = \{n \in \mathbb{N} \, \vert \, M-m_0\leq n \leq N-m_0\}$ for the set of indices of the grid associated with the values of unknown dynamics and $\mathcal{I}_M := \{n \in \mathbb{N} \, \vert \,  M-M_0\leq  n < M-m_0\}$  for the set of indices for supplied initial dynamics. The linear} system \eqref{dyna-eqn} may be written in compact matrix-vector form:
\begin{equation}
B \hat{\fb} = h^{-1} A\xb - \hat\gb,
\label{matrix-eqn}
\end{equation}
{where $A$ is the $N_M \times (N+1)$ 
matrix of coefficients for $\alphab$ corresponding to $\xb_{n-m}$ in \eqref{dyna-eqn};
the matrix $B$ is an $N_M \times N_M$ banded lower-triangular matrix with
its diagonal entries given by $\beta_{m_0}$ and the $k$-th subdiagonal entries given by
 $\beta_{m_0+k}$ for $k=1, 2, ..., M_0-m_0$;
$ \hat{\fb} \in \R^{N_M\times d}$  is the ordered vector of unknowns $\{\hat{\fb}_n\}_{n \in \mathcal{I}}$;}
{and $\hat\gb = (\hat\gb_M, \hat\gb_{M+1}, \ldots, \hat\gb_{N})  \in \R^{N_M\times d}$ is defined as} 
\begin{equation*}
{
\hat\gb_{n} = 
\begin{cases} \displaystyle  
 \sum_{\substack{{m\geq n-M_0} \\ {m\in \mathcal{I}_M } } }
 \beta_{n-m} {\hat \fb}_{m}, \quad & \text{if }\;  n  \in  M_0 + \mathcal{I}_M,\\
0,& \text{if }\;  n \in \mathcal{I} \setminus \{M_0+\mathcal{I}_M\},
\end{cases}
}
\end{equation*}
{which can be generated from the assumed, suitably approximated starting values $\{\hat{\fb}_n\}_{n \in \mathcal{I}_M}$.}
{
We note that since $\beta_{m_0}\neq 0$, $B^{-1}$ always exists so that \eqref{matrix-eqn} is solvable whenever the right hand terms are prescribed.} 
 \subsection{{Connection to Machine Learning-based Data-driven Discovery and LMNet}}\label{sec:connect}
 {To see how the theory developed in this work is connected to the increasingly popular machine learning based data-driven discovery of dynamics, {we} briefly recall the relevant learning problems here. For more extensive works on machine learning, we refer to  \cite{bishop06,goodfellow2016deep,mohri18,murphy12,shalev2014}.}
 
 {
In a generic supervised machine learning setting of learning an unknown function $\fb$, one often
assumes knowledge of $\tilde{N}$ samples of input-output data, $\mathcal{D}= \{(\xb_n, {\fb}(\xb_n))\}_{n=1}^{\tilde{N}}$. This sample dataset is often divided into sets of training and test sets, and one attempts to
find a neural network (NN) representation of $\fb$, say $\fb_{NN}$, through
an empirical loss minimization over the training set. We let $\tilde{\xb}$ and $\tilde{\fb}$ denote an ordered subset of $\tilde{K}\leq\tilde{N}$ data, so that $(\tilde{\xb}_k, \tilde{\fb}_k) = \big(\xb_{n_k}, \fb(\xb_{n_k})\big) \in \mathcal{D}.$
The loss is  a suitably-defined function $\ell(\tilde{\xb}, \tilde{\fb}, \fb_{NN})$ measuring a distance between $\fb(\xb_{n_k})$ and $\fb_{NN}(\xb_{n_k})$ for each $k=1, 2, \ldots, \tilde{K}$.  When evaluated over only training data, this loss leads to the training error.
The desired goal is to learn $\fb_{NN}$ that not only minimizes the loss in
the training set (i.e., the training error), but also achieves a small loss in the remaining test samples  (i.e., the generalization error). }

{
In the setting of dynamics discovery, it is important to note that the dynamics, or output, data is not given directly.
Instead, only the state, or the input,  is provided, and information on the true dynamics $\fb$  is inferred by
constraining the data to conform with some dynamical system. For {the} LMM discretization of the dynamics {given by \eqref{dyna-eqn}, } conformity is achieved by minimizing the {error associated with the LMM system, which we call the LMM residual.} 
A total loss function of the optimization problem may be effectively viewed  as}
$${\cT(\tilde{\xb}, \tilde{\fb}, \fb_{NN})=\tilde{\ell}(\tilde{\xb}, \tilde{\fb})  + \ell(\tilde{\xb},\tilde{\fb} , \fb_{NN} ) }
$$
{where the loss $\tilde{\ell}$ is an increasing function of the {LMM residual} and vanishes at the origin, e.g., 
$\tilde{\ell}(\tilde{\xb}, \tilde{\fb})= \|
B \tilde{\fb} - h^{-1} A\tilde{\xb} + \hat\gb\|_2^2$. 
A network approximation with sufficient accuracy would attempt to 
conform with the discretized LMM dynamics by
minimizing the LMM residual to find the unknown data $\tilde{\fb}$. Alternatively, as done in LMNet, the neural network approximation may be supplied to {the LMM residual} $\tilde{\ell}$, where the initial dynamics in $\hat\gb$ are also learned.
If the approximation can be as accurate as desired, we would be led to the idealized setting }
{that as the network is trained more, given sufficient width, the neural network would converge to $\hat\fb$, where $\hat{\fb}$ satisfies \eqref{matrix-eqn}.}

{
Naturally, due to other practical considerations as well as the finite approximation power of the neural networks, more general loss functions, regularization techniques, and network architectures  may also be taken into account, see Section~\ref{sec:conc} for further discussions. Our main focus here is to illustrate the impact of using different LMM on the learning process by developing a rigorous mathematical theory of consistency, stability and convergence for the dynamics discovery, beginning with the highly idealized setting {of exact state data}.
}

\subsection{Truncation Error and Consistency}\label{sec:trunc-error}

LMMs for discovery inherit the truncation error of solving ordinary differential equations with LMMs. Indeed, truncation error is specific to the discretization of the continuous problem; therefore, the truncation error $\taub_h$ of a scheme for dynamics discovery remains the same as that for solving an ordinary differential equation for the state defined by \eqref{eq:taubhdef}.   
However,  in addition to inheriting the same concept of consistency from Section \ref{sec:prob-intro}, Definition \ref{def:consistency},
we also introduce some strengthened notions of consistency for dynamics discovery. We complement these concepts later on with refined notions of stability for a more nuanced discussion of convergence for discovery using LMMs. Consistency and its strengthened forms are defined below.

\begin{newdef}[Consistency for Dynamics Discovery]\label{def:consistencynew}{  {An LMM} is consistent with the differential equation for dynamics discovery provided $\norm{\taub_h}_\infty \rightarrow 0$ as $h \rightarrow 0$,  and it is strongly consistent if $\norm{\taub_h}_1 \rightarrow 0$ as $h \rightarrow 0$.  Furthermore, a method is consistent of degree $k$, for $k\geq 1$, provided $N^{k-1}\norm{\taub_h}_\infty \rightarrow 0$ as $h \rightarrow 0$}. 
\end{newdef}

\begin{remark}\label{remark-new}
With the Definition \ref{def:consistencynew}, all LMMs having at least $k$-th order truncation error are consistent of degree at least $k$.  Moreover,  since
 \[ \norm{\taub_h}_1 = \sum_{n=M}^{N}  |(\taub_h)_n|
 \leq N \norm{\taub_h}_\infty 
 ,\]
 LMMs having at least {second}-order truncation are automatically  {consistent of degree 2 and thus} strongly consistent.
\end{remark}

Following from the classical truncation error analysis for LMMs, we have the algebraic criteria
for the consistency.

\begin{lemma}[Consistency]\label{claim:consistent}
{
A linear multistep method scheme for dynamics discovery is consistent provided 
that $\rho(1)=0$ and $\rho'(1)=\sigma(1)$. Furthermore, it is consistent of degree $k$ if it is
order ${k}$ in the sense of Definition~\ref{def:err}, that is, 
if $ \rho(e^z)-z \sigma(e^z)=O(z^{{k+1}})$ as $z\to 0$.}
\end{lemma}

\subsection{Stability and the Root Condition for Discovery}\label{sec:stability}
 {In this section we develop stability in a similar spirit as in Section \ref{sec:prob-intro} but also introduce more refined notions of stability for convergence analysis. For discovery, the main distinction from theory for solving the forward problem is that now we consider perturbations to} the recovered dynamics
as opposed to the integrated states for the numerical solution of  the differential equation. To begin
we introduce a linear operator given by
\begin{equation}\label{eqn:res-discovery}
{(\hat{R}_h \hat{\fb})_n := }
{ \sum_{m=m_0}^{M_0} \beta_m \hat{\fb}_{n-m},\quad n=M,M+1,\ldots, N.}
\end{equation}
{Notice $(\hat R_h \hat{\fb})_n$ arises from its forward counterpart \eqref{eqn:res-forward} with the reduced $\betab$ notation.}
\begin{newdef}[
Stability for Dynamics Discovery]\label{def:stable-new2}
A linear $M-$step method for the dynamics discovery is called 
stable on $[ \, a,b\, ]$ provided there exists a constant $K<\infty$, not depending on {$N$}, such that, for any two grid functions $\ub, \vb \in \Gamma_h[ \, a,b\, ]$, 
we have 
\begin{equation*}
\norm{\ub - \vb}_\infty \leq K \left(\max_{{i \in \mathcal{I}_M }}\abs{\boldsymbol{u}_i - \boldsymbol{v}_i} + \norm{\hat R_h (\boldsymbol{u} - \boldsymbol{v})}_\infty  \right).
\end{equation*}
\end{newdef}

\begin{newdef}[{Marginal}
Stability for Dynamics Discovery]\label{def:stable-new}
A linear $M-$step method for the dynamics discovery is called  {marginally}
stable on $[ \, a,b\, ]$ provided that there exists a constant $K<\infty$, not depending on
{$N$}, such that, for any two grid functions $\ub,\vb\in \Gamma_h[ \, a,b\, ]$,
we have 
\begin{equation*}
\norm{\ub - \vb}_{{\infty}} \leq K \left(\max_{{i \in \mathcal{I}_M}}\abs{\boldsymbol{u}_i - \boldsymbol{v}_i} + \norm{\hat R_h (\boldsymbol{u} - \boldsymbol{v})}_1  \right).
\end{equation*}
\end{newdef}

\begin{newdef}[{Weak} Stability of Degree ${-k}$ for Dynamics Discovery]\label{def:stable-new3}
A linear $M-$step method for the dynamics discovery is called {weakly} stable of degree ${-k}$ {for ${k\geq 2}$} on $[ \, a,b\, ]$ provided that there exists a constant $K<\infty$, not depending on {$N$}, such that, for any two grid functions $\ub, \vb \in \Gamma_h[ \, a,b\, ]$, 
we have 
\begin{equation*}
\norm{\ub - \vb}_\infty \leq K   \left(
{N^{{k-2}}}
\max_{{i \in \mathcal{I}_M}}\abs{\boldsymbol{u}_i - \boldsymbol{v}_i} +  {N^{{k-1}}} \norm{\hat R_h (\boldsymbol{u} - \boldsymbol{v})}_\infty  \right).
\end{equation*}
\end{newdef}

{In all cases, the norm on the left-hand-side is taken over the learned 
components $\{\ub_n\}_{{n\in\mathcal{I}}}$ and $\{\vb_n\}_{{n\in\mathcal{I}}}$. This convention is used in the rest of the paper.}
Note that {we choose to use negative degree values so that more negative degrees correspond to weaker stability.}
 Similar to the observation given in Remark \ref{remark-new},  {we see that weak stability of degree $-2$ follows from the marginal stability in Definition \ref{def:stable-new}. 
 }

We would like to turn the property of stability  into an algebraic condition as for the case of numerical solution to ODEs. For the forward problem, the algebraic root condition (Definition \ref{def:root-condition}) serves this purpose; however, for the inverse problem, we require a more subtle treatment of the root condition to capture the nuances in stability for dynamics discovery.

\begin{newdef}[Strong Root Condition \cite{agarwal00,strelitz77,dorf11, atkinson11}
]\label{def:strong-root-condition}
A polynomial satisfies the strong root condition provided the roots of the polynomial have magnitude less than 1.
\end{newdef}

Likewise, we also generalize the above root conditions.
 \begin{newdef}[$k^{th}$-multiplicity Root Condition]\label{def:k-root-condition}
A polynomial satisfies the root condition of degree {$k \in \mathbb{N}$} provided the roots of the polynomial do not exceed magnitude 1, and those of magnitude 1 have multiplicity no larger than $k$. 
\end{newdef}


\begin{remark}
One may view the conventional (algebraic) root condition (Definition \ref{def:root-condition}) and the strong root condition
(Definition \ref{def:k-root-condition}) as special cases of the $k^{th}$-multiplicity root condition of Definition \ref{def:k-root-condition} with $k=1$ and $k=0$ respectively. 
The strong root condition has been used in the numerical analysis, control theory, and linear recurrence relation literature for study of relative stability for LMM
as time integrators and asymptotic properties associated with the linear recurrence relations \cite{agarwal00,strelitz77,dorf11,atkinson11}. 
\end{remark}

Naturally,
we can see that the notions of stability for discovery for LMMs are tied to the bounds on the solutions to the linear recurrence equations determined by the coefficients {$\betab$}.  
We now relate them to the root conditions. Notice that while the stability in
Theorem \ref{thm:root-state} for numerical integration of the given dynamics
 is concerned with the first characteristic polynomial $\rho(r)$,  the stability in
Theorem \ref{thm:stability-new} for the discovery of dynamics is 
concerned with the second characteristic polynomial $\sigma(r)$ { defined by \eqref{eqn:char}. 
More precisely, the root condition can be stated for a reduced second characteristic polynomial
\begin{equation}
\hat{\sigma}(r)=\sum_{m=m_0}^{M_0} \beta_m r^{M_0-m}. 
\label{eq:hatsigma}
\end{equation}}Hence, we see a fundamental difference in the two stability notions. {The dependence of stability on 
$\sigma(r)$  (or $\hat{\sigma}(r)$) might be unexpected as it has not appeared in the numerical differential equation literature.
However, it is also not surprising given the inverse problem nature of using LMMs for dynamics discovery. }

\begin{newthm}[Stability for Discovery]\label{thm:stability-new}
A linear multistep method for discovery of dynamics is 
{stable provided that the second characteristic polynomial $\sigma(r) $
or 
the reduced $\hat{\sigma}(r)$
satisfies the strong root condition in Definition \ref{def:strong-root-condition}.
Likewise,
an LMM for discovery of dynamics is marginally stable provided that $\sigma(r)$ 
or
$\hat{\sigma}(r)$
satisfies the algebraic root condition in Definition \ref{def:root-condition}.
Furthermore,
an LMM for discovery of dynamics is {weakly} stable of degree $-k$ (for $k\geq 2$) provided that
$\sigma(r) $ or  $\hat{\sigma}(r)$  
satisfies 
the $(k-1)^{th}$-multiplicity} root condition in Definition \ref{def:k-root-condition}.
\end{newthm}

\begin{proof} 
{Let $\hat\eb=\ub-\vb$, {where $\ub, \vb \in \Gamma_h[0,T]$ are both generated by solving the LMM \eqref{matrix-eqn}. By setting
$\rb=\hat R_h (\eb)$} with the operator $\hat{R}_h$  defined in \eqref{eqn:res-discovery}, we have
$$
\sum_{m=m_0}^{M_0} \beta_m \hat\eb_{n-m}
=\rb_n ,\quad n=M,M+1,\ldots, N.
$$
By standard recurrence and linear algebra theory \cite{agarwal00,gautschi}, 
 the {difference $\hat\eb$ can be  determined by the companion matrix of the above recurrence relation, denoted by $\mathcal{Z}$. This matrix is an $(M_0-m_0)\times(M_0-m_0)$  matrix with its first row given by {$-\left(\frac{\beta_{m_0+1}}{\beta_{m_0}},\frac{\beta_{m_0+2}}{\beta_{m_0}},\ldots,\frac{\beta_{M_0}}{\beta_{m_0}}\right)$}, and {the} rest of the rows are of the form $(\mathbf{I},\mathbf{0})$ where $\mathbf{I}$ is the identity matrix {of size $(M_0-m_0 -1)\times (M_0-m_0-1)$ and }
$\mathbf{0}$ is the zero column vector in $\R^{M_0-m_0-1}$}. 
 The matrix  $\mathcal{Z}$ is 
associated with a characteristic polynomial given by $\hat\sigma(r)$ that shares the
same set of roots as that of $\sigma(r)$, except a possible root at $0$
.}

{To consider the propagation of the {difference $\hat\eb$, we form the matrix} $\textbf{E}_n=Z\textbf{E}_{n-1}+ \Rbf_n$
where $\textbf{E}_n\in\mathbb{R}^{(M_0-m_0)\times d}$ has its rows given by $\{\hat\eb_{n-k}\}_{0\leq k < M_0-m_0}$, and 
 $\Rbf_n\in\mathbb{R}^{(M_0-m_0)\times d}$ has its first row given by the vector 
$\beta_{m_0}^{-1}\rb_n$, and all}
{subsequent rows by zeros.}
{
Then 
$$\textbf{E}_n=Z^{n-M+1} \textbf{E}_{M-1} + \sum_{k=M}^n Z^{n-k} \Rbf_k$$ where $\textbf{E}_{M-1}$ is given by the initial data
$\{\hat\eb_k\}_{k\in \mathcal{I}_M}$. 
{Thus, stability} is equivalent to
$$\max_{M\leq n\leq N}
\|\textbf{E}_n\|_\infty \leq K(\|\textbf{E}_{M-1}\|_\infty + \max_{M\leq n\leq N} \|\Rbf_n\|_\infty),$$
which is implied  by
\[ \sum_{{n=1}}^{{N_M}}  \norm{\mathcal{Z}^n}_{\infty} \leq K^* < \infty,\]
or equivalently the strong root condition. {Meanwhile, marginal stability}
is equivalent to
$$ \max_{M\leq n\leq N} \|\textbf{E}_n\|_\infty \leq K(\|\textbf{E}_{M-1}\|_\infty + \sum_{M\leq n\leq N} \|\Rbf_n\|_\infty),$$
 which is implied  by
\[ \max_{{1 \leq n\leq N_M}} \norm{\mathcal{Z}^n}_{\infty} \leq K^* < \infty.\]
We thus only need the algebraic root condition. }

{Likewise, we can argue that {weak stability of degree $-k$ is implied by}
\[ \max_{{1\leq n\leq N_M}} \norm{\mathcal{Z}^n}_{{\infty}} \leq K^*  N^{{k-2}} < \infty,
\quad\text{and}\quad \sum_{{n=1}}^{{N_M}}  \norm{\mathcal{Z}^n}_{\infty} \leq K^* N^{{k-1}},
\]
which is equivalent to the {${(k-1)}^{th}$ multiplicity} root condition.} 
\end{proof}


\subsection{Error Analysis and Convergence}
{In this section, we use the truncation error to study the error for discovery, including defining convergence and the order of approximation of LMM schemes for discovery.}

\begin{newdef}[Convergence {and Order of Approximation} for Discovery]\label{def:convergence-disc}
Consider the initial value problem \eqref{model0} 
{discretized by an $M-$step LMM given by \eqref{dyna-eqn}.} {Let $\fb, \hat{\fb} \in \Gamma_h[ \, a,b\, ]$, where $\fb$ is the exact grid function on the $N+1$ grid points 
 $\{\fb_n=f(x(t_n))\}$  and $\hat{\fb}$ the approximation solved from \eqref{dyna-eqn}. {The LMM is convergent {for dynamics discovery} if}
 $
 \norm{\fb-\hat{\fb}}_\infty \rightarrow 0$ as $h\rightarrow 0$ whenever $\max_{ {i\in \mathcal{I}_M} }\abs{\fb_i-\hat{\fb}_i}  \to 0.$
Moreover, if  $\norm{\fb-\hat{\fb}}_\infty = c h^p+O(h^{p+1})$ for some constant $c$, then $p$ is called the convergence order, or alternatively, the order of approximation for dynamics discovery.} 
\end{newdef}


Using the introduced notions of consistency and stability, we now present convergence theorems for dynamics discovery.

\begin{newthm}[Convergence Theorems for Discovery {I}]\label{thm:conv-disc}
{Consider the dynamical system \eqref{model0} discretized by an $M-$step LMM given by 
\eqref{dyna-eqn}.} {Let $\fb, \hat{\fb} \in \Gamma_h[ \, a,b\, ]$, where $\fb$ is the exact grid function on the $N+1$ grid points 
 $\{\fb_n=f(x(t_n))\}$  and $\hat{\fb}$ the approximation solved from \eqref{dyna-eqn}.}
  Then,
\begin{equation}\label{eqn:err}
B( \hat{\fb} -\fb)= \taub_h + \boldsymbol{g}_h ,
\end{equation}
where {$\taub_h$ is the local truncation error of the scheme, $\gb_h = (\gb_M, \gb_{M+1}, \ldots, \gb_{N})$} is given by
\begin{equation*}
{ 
(\boldsymbol{g}_h)_n =
\begin{cases}
\displaystyle
\sum_{\substack{{m\geq n-M_0} \\ {m\in \mathcal{I}_M } } }
 \beta_{n-m}  ( {\fb}_{m} - {\hat \fb}_{m}),
\quad & \text{if }\;   n  \in  M_0 + \mathcal{I}_M,\\
0,& \text{if }\;  n \in \mathcal{I} \setminus \{M_0+\mathcal{I}_M\},
\end{cases}
}
\end{equation*}
{Moreover,  in the senses of consistency outlined in Definition \ref{def:consistencynew}
 and stability in Definitions  \ref{def:stable-new2}-\ref{def:stable-new3}, 
if an LMM is consistent and stable, or strongly consistent and marginally stable, then it is convergent for dynamics discovery {in the sense of} Definition~\ref{def:convergence-disc}.
Furthermore, if it is consistent of degree $k$ and weakly stable of degree $-k$,
 then provided that 
  $N^{k-2}\max_{ {i\in \mathcal{I}_M} }\abs{\fb_i-\hat{\fb}_i}  \to 0$  as $h\rightarrow 0$, we also have
 $\norm{\fb-\hat{\fb}}_\infty \rightarrow 0$.
}
\end{newthm}

\begin{proof} 
{By Equation \eqref{dyna-eqn} and the truncation error defined in Equation \eqref{eq:taubhdef}, we have}
\begin{eqnarray*}
 {\sum_{m=0}^M  \frac{1}{h}\alpha_m \xb_{n-m} -  \sum_{m=m_0}^M  \beta_m \hat\fb_{n-m} }
 &=& \boldsymbol{0},\\
{\sum_{m=0}^M \frac{1}{h}\alpha_m \xb_{n-m} -   
\sum_{m=m_0}^M  \beta_m \fb_{n-m}} &=&  (\taub_h)_{n}.
\end{eqnarray*}
{Subtracting the equations,  we observe}
\begin{equation}\label{eqn:errorrecur}
{
\sum_{m=m_0}^M \beta_m 
( \hat\fb_{n-m} - \fb_{n-m} )   =  (\taub_h)_n, \quad n=M,M+1,\ldots, N,
}
\end{equation}
{or equivalently
$B (\fb - \hat{\fb})= \taub_h +\boldsymbol{g}_h$
where $\hat{\fb} - \fb$ on the left side refers to those components indexed in $\mathcal{I}$.}

{Now,  
by the definitions of stability given in Definitions \ref{def:stable-new2} and \ref{def:stable-new},}
there exists a constant $K_W<\infty$ independent of $h$,  for $h$ sufficiently small, such that
\[
\begin{aligned}
\norm{ \fb - \hat{\fb} }_{\infty} 
 \leq  K_W \left(
 { \max_{{i\in \mathcal{I}_M}}\abs{\fb_i - \hat{\fb}_i}
 +  \norm{ 
 \taub_h }_W }
  \right),
\end{aligned}
\]
{where $W=\infty$ or $W=1$, if the LMM is stable or marginally stable, respectively.
Thus,  by Definition 
\ref{def:consistencynew} on consistency and strong consistency,  we have
 $\norm{\fb - \hat{\fb}}_\infty \to 0$
as $h\rightarrow 0$. }
{Likewise,
 if the LMM is 
{weakly} stable of degree $-k$,
 then
\[
\begin{aligned}
\norm{ \fb - \hat{\fb} }_{\infty} 
& \leq K  \left(  N^{k-2}
{\max_{{i\in \mathcal{I}_M}}\abs{\fb_i - \hat{\fb}_i}
 + N^{k-1}  \norm{
 \taub_h }_\infty }
\right).
\end{aligned}
\]
By the definition on the consistency of degree $k$,  
 together with the assumption on the initial data that
  $N^{{k-2}} \max_{{i \in \mathcal{I}_M}}\abs{\fb_i-\hat{\fb}_i} \to 0$,
convergence also follows.}
\end{proof}



{Theorem 
\ref{thm:conv-disc}
 states that for LMM based dynamics discovery, convergence follows from both stability and consistency,
 as in the case of LMM-based time integration}. 
 The equation \eqref{eqn:err} shows the interplay between {the stability aspect of } solving the system, manifested in $B^{-1},$ and the {consistency component of} truncation error, $\taub_h$, from discretization of the differential equation.

{We note that Theorem~\ref{thm:conv-disc} contains only
sufficient conditions for convergence.
There are examples of LMMs that are consistent and marginally stable, but not strongly consistent {nor} stable, which may  still  be convergent for dynamics discovery. An example is the LMM with $\rho(r)=r^2-1$ and $\sigma(r)=(r+1)/2${; convergence for this LMM} can be checked using calculations similar to that presented in the proof of Corollary~\ref{cor:conv}.  Nevertheless, in the spirit of the {Dahlquist Equivalence Theorem}, we also have the following result establishing consistency and some form of stability from convergence.}

\begin{newthm}[Convergence Theorems for Discovery {II}]\label{thm:conv-disc-new}
{Consider the dynamical system \eqref{model0} discretized by an $M-$step LMM given by 
\eqref{dyna-eqn}. If the LMM is convergent for dynamics discovery
{in the sense of} Definition~\ref{def:convergence-disc}, {then it is consistent and marginally stable in the senses of Definitions~\ref{def:consistencynew}
and \ref{def:stable-new}.}
}
\end{newthm}

\begin{proof} 
{
The proof is similar to {its} classical counterpart. {Consider the} special {cases of ODE}
$\frac{d}{dt}\xb(t) = 0$ and $\frac{d}{dt}\xb(t) = 1$, respectively with $\xb(a)=0$. {If the LMM is convergent in the sense of Definition \ref{def:convergence-disc}, then} the dynamical system
\eqref{dyna-eqn} leads to {a} linear recurrence relation with constants $\rho(1)$ or $\rho'(1)$, respectively, serving as inhomogeneous {terms on the} right hand side. Since the LMM is convergent, the learned dynamics approach $0$ or $1$, respectively{. Thus, as} $h \rightarrow 0,$  we get $\rho(1)=0$ and $\rho'(1)=\sigma(1)$. Consequentially, the consistency conditions from Lemma~\ref{claim:consistent} are satisfied. Using the theory on the linear recurrence relations given in the proof of Theorem \ref{thm:conv-disc}, 
{in order for}  $ \norm{\fb-\hat{\fb}}_\infty \rightarrow 0$ as $h\rightarrow 0$ whenever $\max_{ {i\in \mathcal{I}_M} }\abs{\fb_i-\hat{\fb}_i}  \to 0$, 
there must exist some constant $0 < K < \infty$ such that $\max_{1\leq n \leq N_M} \|\mathcal{Z}^n\| < K$ as $N_M \nearrow \infty$. For this bound to exist, the root condition must be satisfied, and hence the method must be marginal stable.}
\end{proof}

{As seen from the proof of Theorem \ref{thm:conv-disc}, under some assumptions on the initial
dynamics, we immediately get the order of convergence for LMM-based dynamics discovery.}

\begin{newthm}[Order of Convergence]\label{thm:order}
{
Given an LMM with a truncation error of order $k$ with $k\geq 1$, i.e.,  $\norm{\taub_h}_\infty= O(h ^{k})$,  as in Definition~\ref{def:err}. 
Then, as $h\to 0$,  we have
 $\norm{\fb-\hat{\fb}}_\infty = O(h^{k})$
if the LMM is 
stable and
$\displaystyle \max_{{i\in \mathcal{I}_M}}
 \abs{\fb_i - \hat{\fb}_i} = O(h^{k})$. Moreover, 
  provided that
$\displaystyle \max_{{i\in \mathcal{I}_M}}
 \abs{\fb_i - \hat{\fb}_i} = O(h^{k-1})$, we have
 $\norm{\fb-\hat{\fb}}_\infty \leq C h^{k-1}$
 if it is {marginally} stable or 
  $\norm{\fb-\hat{\fb}}_\infty \leq C h^{k+1-s}$
  if it is weakly stable of degree $-s$ with $k\geq \max\{s,2\}$.
 } 
\end{newthm}

\begin{remark}\label{rmk:conv2}
{
The different notions of stability affect the order of convergence for dynamics discovery. These refinements motivate accompanying definitions for the degree of consistency in Definition~\ref{def:consistencynew}, whereas traditionally the {order of error matches with the} order of convergence (see Definition~\ref{def:err}).
{For dynamics discovery,} the order of convergence and degree of consistency {matches} for {strongly} stable schemes.}
 {While this might not hold generically for marginally or weakly stable LMMs, resulting in possible lower order of convergence than the degree of consistency. {We} show later that, for some cases such as AM-1, the same order can still be maintained.}
 \end{remark}


\section{Application to AB, AM, and BDF}\label{sec:app}

We now apply the general theorem on LMM for the dynamics discovery to three popular special classes of methods-- Adams-Bashforth (AB), Adams-Moulton (AM), and Backwards Differentiation Formula (BDF).

\subsection{Consistency of AB, AM and BDF}\label{sec:consis}

It is {well-known} that 
{the Adams Family schemes and BDF are consistent as time integrators. Specifically, AB-$M$ and BDF-$M$ have order of error $M$, while AM-$M$ has order of error $M+1.$} As a result, {these} three classes of LMM methods 
  remain  consistent for dynamics discovery. {Moreover, as a consequence of the order of error, AB-$M$ and BDF-$M$ are consistent of degree $M$, and the AM-$M$ schemes are consistent of degree $M+1,$} as noted in Remark \ref{remark-new}. Indeed, the latter fact is crucial to the convergence of AM-1.

\begin{newthm}[Consistency of AB, AM and BDF for Dynamics Discovery]\label{thm:consist-dynamics}
The 
{AB-$M$, AM-$M$ and BDF-$M$ schemes}
 are all consistent for dynamics discovery. Furthermore, 
{AM-1 is consistent of degree $2$ and thus} strongly consistent .
\end{newthm}

\subsection{Stability and Convergence of AB, AM and BDF}
\label{sec:stability-methods}
\begin{newthm}\label{thm:scheme-stability}
With the notions of stability defined in Definitions \ref{def:stable-new} and \ref{def:stable-new2},
\begin{enumerate}
\item {BDF-$M$ for all $M \geq 1$, AB-$M$ for $1 \leq M \leq 6$, and AM-0}  are 
 stable;
\item {AM-1} is 
{marginally} stable and thus {weakly stable of degree $-2$};
\item {AB-$M$} for $7 \leq M \leq 10$ and {AM-$M$ for $M \geq 2$} are unstable.
\end{enumerate}
\end{newthm}

The proof of Theorem \ref{thm:scheme-stability} is given in Section \ref{sec:char-poly}.

\begin{cor}\label{cor:conv}
{BDF-$M$  for all $M\geq 1$
are convergent, with
convergence order $M$.}
{AB-$M$ for $1 \leq M \leq 6$ are convergent, with
convergence order $M$.
 AM-0  is  convergent with first-order convergence. 
AM-1 is convergent with second-order convergence if we have second order error
on the initial data}. 
\end{cor}
\begin{proof}
{The conclusions of Corollary \ref{cor:conv} on the convergence of LMM schemes under consideration follow immediately from the application of 
Theorems \ref{thm:consist-dynamics}, \ref{thm:scheme-stability}, and \ref{thm:order}.}
 {The order of convergence follows, with the exception of 
 {AM-1}. Indeed, a direct application would imply only first order convergence due to its degree-1 marginal stability. However, we note in this special case, 
 the recurrence relation \eqref{eqn:errorrecur}
  is given by $
\hat \fb_n- \fb_n=- (\hat\fb_{n-1} - \fb_{n-1})+2(\taub_h)_n$.
 Using the error expansion given in Definition~\ref{def:err},
 the leading order of  $\hat \fb_n- \fb_n$
 of the form
 \begin{align*}
 &   (-1)^{n-j} (\hat\fb_0-\fb_0 ) h^2 + \sum_{j=1}^{n}  (-1)^{n-j}  \eb(t_j) h^2    \\
 &\quad \approx O(h^2) + 
 \left\{
 \begin{array}{ll}
 \displaystyle
 \sum_{j=1}^{k} \eb'(t_{2j}) h^{3},  &\; \text{ if $n=2k$,}\\
  \displaystyle
 \sum_{j=1}^{k} \eb'(t_{2j+1}) h^{3}  +(-1)^{n-1} \eb(t_1)h^2, &\; \text{ if $n=2k+1$}
 \end{array} 
 \right\}
 \approx O(h^2),
 \end{align*}
 where $\eb(t)=x^{(p+1)}(t)$ is  assumed to be a  smooth function depending on the solution of the exact dynamic $x=x(t)$. 
 Therefore, {given second-order error in the initial data, {AM-1} has second-order convergence  even though it is not a strongly stable method.}
 }
\end{proof}

\begin{remark}
The finite range of instability with respect to the order $M$ for the AB scheme is due to limitation of explicit calculations. We conjecture that the scheme is unstable for all $M \geq 7$.  Interestingly, that $M=6$ is a threshold for stability of the polynomial echoes the stability criterion for the forward problem BDF \cite{henrici62}, for which $M=6$ is also the largest known order method that is stable. Explicit numerical calculation or Routh Arrays (see \cite{dorf11}) are used to show this fact \cite{henrici62, cryer72,fredebeul98}. Schur polynomials have since been used \cite{creedon75} to show a generalized stability argument for $M\geq 13$ \cite{fredebeul98}. We leave open a generalized stability result for $M\geq7$ using the polynomial roots, but we have validated numerically the instability for $ 7 \leq M\leq 20$.
\end{remark}

\subsection{Verification of Root Conditions for AB, AM and BDF}\label{sec:char-poly}

We now verify, for the three classes of LMMs, the root condition holds for cases stated in Theorem \ref{thm:scheme-stability}.

We begin by calculating 
the roots of the second characteristic
polynomial associated with AB and AM since {$\sigma(r)=\hat{\sigma}(r)$ in both cases.
We first present some results for AB-$M$ and AM-$M$ with $1 \leq M \leq 10$ as 
computational evidence (with exact symbolic computation).
 We have also numerically validated instability for  AB-$M$ for $11 \leq M \leq 20$ as well and expect instability to persist for all $M \geq 7$.  However, there is no theoretical proof so far.
  For AM-$M$, a general instability result for $M\geq 2$ is proved in Lemma \ref{claim:exp-growth}.} 

 Fix $M\in \mathbb{N}$ and  $\tilde\Lambda \in \{\Lambda_0, \Lambda_1\}$, where we recall from Section \ref{sec:prob-intro} that
$\Lambda_0 = \{-M, \ldots, 0\}$   and  {$\Lambda_1 = \{-M,  -M+1,\ldots, -1\}$}.
Exchanging the integral and the summand in the formula for the Lagrange interpolating polynomial, one can observe that finding the roots of the second characteristic polynomial is equivalent to choosing $r \in \mathbb{C}$ satisfying a mean-zero equation. {That is, for  $\ell_x^h(u; \tilde\Lambda)$
defined in \eqref{eqn:lag-basis-ele-uniform}, we have}
\begin{equation}
\sum_{x \in \tilde\Lambda}\int_0^1 \ {\ell_x^h}(u; \tilde\Lambda)r^xdu  \iff \int_0^1 \sum_{x \in \tilde\Lambda} {\ell_x^h}(u;\tilde\Lambda)r^xdu = 0.\label{eq:roots}
\end{equation}
\vspace{-3mm}
 \begin{table}[H]\caption{Largest Magnitude Roots}\label{tab:roots}\centering
\begin{tabular}{ | c | c | c | c  | c | c | }\hline
Step $M$   & 1 & 2 & 3 & 4 & 5  \\ \hline
AB  & --  & 0.3333 & 0.4663 & 0.6338 & 0.8075   \\
AM & 1.0000& 1.7165 & 2.3658 & 2.9775 & 3.5639  \\\hline\hline
Step   & 6 &7 & 8 & 9 & 10 \\ \hline
AB   & 0.9829& 1.1587 & 1.3345 & 1.5100 & 1.6852 \\
AM & 4.1312  & 4.6851 & 5.2267 & 5.7586 & 6.2820 \\\hline
\end{tabular}
\end{table}

As we see in Table \ref{tab:roots}, which is computed symbolically by Mathematica, 
 the profile of the roots of the characteristic polynomial associated with the different schemes varies significantly. {Equation \eqref{eq:roots} and the data in Table \ref{tab:roots} immediately lead to the following lemma.}
 
\begin{lemma}\label{claim:poly}
Fix $\tilde\Lambda \in \{\Lambda_0, \Lambda_1\},$ and let ${\ell_x^h}(u;\tilde\Lambda)$ be the Lagrange interpolating polynomial defined { in \eqref{eqn:lag-basis-ele-uniform}.}
Then, we can characterize the roots $r \in \mathbb{C}$ of {the second characteristic polynomial as the solution to the equation}
\begin{equation}\label{mean-eqn}
{\int_0^1 \sum_{x \in \tilde\Lambda} {\ell_x^h}(u; \tilde\Lambda)r^xdu = 0.}
\end{equation}
{Moreover, for AM-$M$ with $\tilde\Lambda = \Lambda_0,$ we have
\begin{enumerate}
\item $M=1,$ then the single root satisfies $\abs{r}=1$.
\item $2 \leq M \leq 10$, there exists at least one root $r$ with $\abs{r} > 1.$
\end{enumerate}
and for AB-$M$ with $\tilde\Lambda = \Lambda_1,$ we have
\begin{enumerate}
\item  $1 \leq M \leq 6,$ then $\abs{r}<1$.
\item $7 \leq M \leq 10$, there exists at least one root  $r$ with $\abs{r} > 1.$
\end{enumerate}
}
\end{lemma}


{Let us state some useful properties of the second characteristic polynomial $\sigma(r)$ associated with the AM methods and the corresponding  coefficients of its $B$ matrix.}

\begin{lemma}\label{lemma:am-coeff-results}
For  $M\geq 2,$ the coefficients $\{\beta_m\}_{0}^M$ of the AM method have the properties:
\begin{enumerate}
\item $\beta_1 > \beta_0{>0}$,  \label{beta01}
\item  $\sign(\beta_{m+1}) = -\sign(\beta_m), \ 1 \leq m \leq {M-1}$, and \label{am-sign}
\item { $\beta_0 > |\beta_M|$.}\label{am-extra}
\end{enumerate}
\end{lemma}

\begin{proof} 
Fix $M \in \mathbb{N}$ with $M \geq 2.$ The $M-$step AM scheme has coefficients
\begin{equation}
\beta_{m} =
 \frac{(-1)^m}{m!(M-m)!} \int_0^1 \prod_{\substack{i=0 \\ i \not=m}}^M (u+i-1) \ du,\label{am-coeff}
\end{equation}
for  $ m = 0, 1, \ldots, M.$ 
The coefficients $\beta_0$ and $\beta_1$ are given by
\begin{equation*}
\begin{split}
\beta_0 = \frac{1}{M!} \int_{0}^1 \prod_{i=0}^{M-1} (u+i) du \ \ \text{ and } \ \
\beta_1 = \frac{1}{(M-1)!}\int_{0}^1 (1-u)\prod_{i=1}^{M-1} (u+i) du.
\end{split}
\end{equation*}
{Certainly,} $\beta_0>0$.
Notice 
\begin{align}
\beta_1 > \beta_0 \iff \frac{M}{M+1} \int_0^1 \prod_{i=1}^{M-1} (u+i) du>\int_0^1 \prod_{i=0}^{M-1} (u+i) du.\label{eqn:inductive-hyp-am}
\end{align} 
We prove \eqref{eqn:inductive-hyp-am} by induction.
\
As the base case, $M=2$. For $M=2$, we have $\beta_1=8/12>\beta_0= 5/12$. 
Now assume \eqref{eqn:inductive-hyp-am} holds up to some arbitrary $M \in \mathbb{N},$ with $M > 2$. We will show the result for $M+1.$
\begin{align}
\frac{M+1}{M+2} \int_0^1 \prod_{i=1}^M (u+i) du &= \frac{M+1}{M+2}{\int_0^1} \left( u \prod_{i=1}^{M-1} (u+i ) + M\prod_{i=1}^{M-1}(u+i)\right){du}\\\
&\substack{\eqref{eqn:inductive-hyp-am}\\>} \frac{M+1}{M+2} \left(\int_0^1 \prod_{i=0}^{M-1} (u+i) + \frac{M(M+1)}{M}\prod_{i=0}^{M-1}(u+i)du\right)\label{eqn:induct-use}\\
&= \frac{(M+1) (M+2)}{(M+2)} \int_0^1 \prod_{i=0}^{M-1} (u+i)du\\
&> (M+1)\int_0^1 \frac{u+M}{M+1}\prod_{i=0}^{M-1} (u+i)du= \int_0^1\prod_{i=0}^{M} (u+i)du,
\end{align} 
as desired. Note we used the inductive hypothesis on the second term in \eqref{eqn:induct-use}.  The proof by induction showing for $M\geq 2$, $\beta_1 > \beta_0$ is complete.
To prove  Part \ref{am-sign}, note that the relation of signs between coefficients follows from the sign of the Lagrange basis polynomials in the integrand of the coefficients. For $m \in \{2, 3, \ldots, M\},$ the integrand of \eqref{am-coeff} are of the same sign, and therefore the sign of $\beta_m$ depends only on the multiplier $(-1)^m.$ Hence Part \ref{am-sign} of Lemma \ref{lemma:am-coeff-results} follows.

{Finally, for Part \ref{am-extra}, we note that
\begin{eqnarray*}
|\beta_{M}| & = &
 \frac{1}{M!} \int_0^1  (1-u) \prod_{i=0}^{M-2} (u+i) \ du <   \frac{1}{M!} \int_0^1  \prod_{i=0}^{M-2} (u+i) \ du\\
 & < &  \frac{1}{M!} \int_0^1 (u+M-1)  \prod_{i=0}^{M-2} (u+i) \ du =\beta_0.
\end{eqnarray*}
This completes the proof.}
\end{proof}

{
\begin{lemma}[General Instability of AM $M \geq 2$]\label{claim:exp-growth}
The linear multistep method formed by the Adams-Moulton scheme for $M\geq 2$ does not satisfy the root condition.
\end{lemma}
}

\begin{proof} {
Fix $M\geq 2$ and consider the second characteristic polynomial associated with the Adams-Moulton scheme. We write it as $\sigma(r)=\sum \beta_m r^{M-m}$. From Lemma \ref{lemma:am-coeff-results}, $\beta_1/ \beta_0>1$. Moreover, by construction of the AM method,  $(-1)^{m}\beta_m<0$ for $m\geq 2$. 
}

{For $r>0$ sufficiently large.
\begin{align}
(-1)^M \sigma(-r) &=  (-1)^{2M} r^M \left[ \beta_0  - \beta_1/r +  \sum_{m=2}^M (-1)^m \beta_m r^{-m}  \right].
\end{align}
 Taking the limit as $r \rightarrow +\infty$, we see that $(-1)^M \sigma(-\infty) = \infty$ since $\beta_0 > 0.$ Meanwhile, 
\[(-1)^M\sigma(-\beta_1/\beta_0) =\sum_{m\geq 2}  (-1)^{-m} \beta_m (\beta_1/\beta_0)^{{M-m}} < 0.
\] 
Hence, it follows from the Intermediate Value Theorem that there is at least one  root of $\sigma(r)$ that is real in
 $(-\infty, -\beta_1/\beta_0)\subset (-\infty, -1)$,  violating the root condition.
 The result thus follows.
}
\end{proof}

\begin{newthm}[Root Condition of AB, AM, BDF]\label{thm:root-condition}
{The {strong} root condition for discovery is satisfied by BDF-$M$  for all $M \in \mathbb{N}$, AB-$M$ scheme for $ 1 \leq M \leq 6$, and AM-$M$ for $M=0$.
The algebraic root condition, or the $k^{th}$ root condition with $k=1$, is satisfied  for  {AM-$M$} with $M=1$.
On the other hand, the root condition is not satisfied for  the  {AB-$M$} scheme with $7 \leq M \leq 10$  or the  {AM-$M$} scheme with $M \geq 2$.}
\end{newthm}

\begin{proof} 
{ The case of AM-$0$ is trivial.  Lemma \ref{claim:poly} implies} the results of Theorem \ref{thm:root-condition} for AB-$M$ with $1 \leq M \leq 10$ and for AM-$M$ with $1 \leq M \leq 10$. 
{Furthermore, by Lemma \ref{claim:exp-growth}, the AM-$M$ scheme  for $M \geq 2$ violates the root condition and hence is unstable.
Finally,  BDF-$M$ has $\sigma(r) = r^{M-1}$ and $\hat{\sigma}(r)=1$, for all $M\geq 1$.
 Hence, the root condition is always satisfied for the BDF scheme for arbitrary $M\geq 1$. As a result,  AM-$0$, identical to BDF-$1$, satisfies the root condition as well. }
\end{proof}

Finally, Theorem \ref{thm:scheme-stability} follows directly from Theorems \ref{thm:root-condition} and \ref{thm:stability-new}.

\subsection{{Discussions on the Effect of Initial Conditions}}\label{sec:disc}
{
The theory developed so far is under the assumption that some initial data of the dynamics are provided, which leads to learning the approximated dynamics at later times. One may consider a situation where the some terminal data are given instead. In such cases, the approximate dynamics would be solved backward in time, yielding a
modified system of equations. It is not hard to check that the stability would become dependent on a modified second characteristic polynomial whose roots are the reciprocals of those of $\hat{\sigma}$. Naturally, it is of interest to check root conditions for the
three classes of LMMs as well. For BDF, we clearly see the strong root condition holds as  $\hat{\sigma}(r)=1$. 
For AM-0 and AB-1, the same also hold. Likewise for AM-1, the root condition but not the strong root condition remains true.
For AM-$M$ with $M\geq2$, 
Part \ref{am-extra} of Lemma~\ref{lemma:am-coeff-results} implies the product of the roots of $\hat{\sigma}(r)=\sigma(r)$ is less than one. Therefore, there might be at least one root of the modified second characteristic polynomial outside the unit disc, and hence  instability for these methods is again expected. Interestingly, unlike in the case with initial data where there is not yet rigorous theory but only computational results for the AB methods, one can prove rigorously  in the next lemma a result of instability for the backwards-in-time AB-$M$, $M\geq 2,$ via a similar argument as Part \ref{am-extra} of Lemma~\ref{lemma:am-coeff-results}. 
}

\begin{lemma}\label{lemma:ab-coeff-results}
{
For  $M\geq 2,$ the coefficients $\{\beta_m\}_{0}^{M}$ of the AB-$M$ method satisfy 
$ \beta_0=0$,   and $\beta_1 > |\beta_{M}|$.
}
\end{lemma}

\begin{proof} 
{
$\beta_0=0$ is true by construction. 
For $M \geq {2}$, we have 
\begin{equation}
\beta_{m} =
 \frac{(-1)^{m+1}}{(m-1)!(M-m)!} \int_0^1 \prod_{\substack{i=1 \\ i \not=m}}^{M} (u+i-1) \ du,\label{ab-coeff}
\end{equation}
for $m = 1, 2, \ldots, M+1.$ 
The coefficients $\beta_1$ and $\beta_{M}$ satisfy 
\begin{equation*}
\begin{split}
|\beta_{M} |= \frac{1}{{(M-1)!}}\int_{0}^1 \prod_{i=1}^{M-1} (u+i-1) du
<\frac{1}{(M-1)!} \int_{0}^1 \prod_{i=2}^{M} (u+i-1) du =
\beta_1
\end{split}
\end{equation*}
which completes the proof of the lemma.}
\end{proof}

{From the above, we see that root conditions do not hold for the modified second characteristic polynomial associated with { AB-$M$ with $M\geq 2$}, so that instability would occur when terminal data are supplied.
 In practice, it is often the case that such initial dynamics are represented by neural networks as part of the unknown as well. Thus, the stability in such cases is worthy of further investigation, particularly in conjunction with the approximation properties of the neural networks to be employed. Clearly, the successful runs using neural networks in Figure~\ref{fig:ab-am-bdf0} have good correspondence with those schemes enjoying some stability properties in one or both types of initial/terminal data .
}

\section{Long Time Dynamics Discovery}\label{sec:long-time}
{
In this section, we consider the problem of discovering dynamics of \eqref{model0}
over a variable interval $(0,T),$  with terminal time $1\ll T\to \infty$, and a fixed mesh $h$. Notice by increasing $T$  we increase the number of grid points $N=T/h$; hence we hope to relate our previous studies with variable mesh and fixed domain to this setting.
 For the numerical analysis of time integration,
this study is reminiscent to that of asymptotic stability, which is often treated via the study of linear dynamics \cite{gautschi, suli03, atkinson11}.
}

{
By rescaling time,
$\tilde{t}=t/T$, where $0 \leq \tilde{t} \leq 1,$ and defining $\tilde{x}(\tilde{t})={{x}}(t)$, we have via change of variables that the scaled dynamics $\tilde{{f}}$ may be related to that of the original variables by
\[
\frac{d}{d\tilde{t}}\tilde{{x}}(\tilde{t}) =  T \frac{d}{dt}{{x}}({t}) = T{f}({x}(t)) = T{f}(\tilde{{x}}(\tilde{t})).
\]
Then, if we define $\tilde{{f}}(\tilde{{x}}(\tilde{t}))  = T{f}(\tilde{{x}}(\tilde{t}))=T{f}({x}(t))$,
 the rescaled differential equation becomes
\begin{equation}\label{model0n1}
\frac{d}{d\tilde{t}}\tilde{{x}}(\tilde{t}) = \tilde{{f}}(\tilde{{x}}(t)) ,\ \ 0 \leq \tilde{t} \leq 1 , \ \ \tilde{{x}} (0) = {x}(0)={x}_0 .
\end{equation}
}
{
Now, consider applying the LMM scheme to $\tilde{ {x}}$ using the transformed model problem \eqref{model0n1} with a step size $\tilde{h}=1/N$.  
 Under this rescaling of time, one can check directly the leading 
{truncation error term of an LMM} of order $p$ in the sense of  Definitions \ref{def:trunc} and \ref{def:err} is  
\begin{equation}\label{lt0}C_{p+1}
\tilde{h}^p  \frac{d^{p+1}}{d\tilde{t}^{p+1}}\tilde{ {x}}(\tilde{t})
= C_{p+1} 
\tilde{h}^p  T^{p+1}  \frac{d^{p+1}}{dt^{p+1}}  {x}(t) 
=C_{p+1} T
{h}^p  \frac{d^{p+1}}{dt^{p+1}} {x}(t) .
\end{equation}
}

{
In light of \eqref{lt0}, we can see that the truncation error of the discovered dynamics of \eqref{model0} in the original time scale is a multiple of the truncation error of the rescaled model \eqref{model0n1} by the factor $T^{-1}$. 
Meanwhile, from the analysis of Section \ref{sec:stability}, the error from stability is only directly dependent on  $\sigma(r)$ and $N$, not the specific time domain.
}

{
Using these observations of the effects on consistency and stability, we can deduce the behavior of an LMM in the long-time regime. For a strongly stable $p^{th}-$order LMM,  
the global error behaves like $O\left(T^{-1}T h^p\right)=O(h^p)$ provided that $\max_{t\in(0,T)} |{x}^{(p+1)}(t)|$ remains uniformly bounded as $T$ increases. Hence, we may view strongly stable LMMs as A-stable, in the case of dynamics discovery, for fixed $h$ as $T\to \infty$. This can be seen as another difference with the case of the forward problem of time integration, where the order of A-stable LMMs is known to be limited by 2 due to the celebrated Dahlquist barrier theorems \cite{dahlquist63,gautschi, suli03, atkinson11}.
On the other hand, for unstable methods, the exponential growth in $N$ of the inverse matrix $B^{-1}$ dominates over any gain in accuracy from consistency. Thus, lack of stability leads to 
an exponentially increasing error as $T$ grows linearly. 

{
As a special example, the marginally stable AM-1 
is not  
 stable  for dynamics discovery, but as stated in the Corollary~\ref{cor:conv} and the derivation in its  proof, we can use the rescaling to get the global error  in the form $O(T^{-1}  T h^2)=O(h^2)$,  Thus, we expect  AM-1,  for a fixed $h$,  to have a constant error as $T$ increases, which is supported by numerical experiments presented in the next section.
}

{
To recap, from the analysis in this section, for dynamics discovery,  {BDFs enjoy asymptotic  stability for  a fixed time step size $h$ as $T$ increases. Same holds for AB-$M$, at least for a small value of $M$ that enjoys the stability as $h\to 0$ for a given terminal time. While this also holds for AM-1, it does not hold for  AM-$M$ with $M\geq 2$.}  As shown in Figure \ref{fig:long-time}, the errors from AB and BDF remain fixed across various values of $T$, while the AM methods yield exponential growth of error in $T$ for $M\geq 2$. 
}

\section{Numerical Experiments}\label{sec:numerics} In this section, we {provide numerical solutions to the linear systems associated with} each of the studied multistep methods and show numerical evidence consistent with the theoretical findings. We limit ourselves to {the idealized setting of numerically exact states considered for the theoretical analysis and to} low dimensional dynamic systems for the sake of illustration and benchmarking. {In addition, we also take the initial data for the dynamics to be exact.} For a model problem, we consider the 2D Cubic System, a nonlinearly damped oscillator, specified as in \cite{mnn, brunton16}. 
\begin{equation} \label{model1}
\begin{cases}
& \dot{x}_1= - 0.1 \ x_1^3 + 2.0 \ x_2^3,\\
& \dot{x}_2= - 2.0 \ x_1^3 -0.1 \ x_2^3,\\
&  {x_1(0) = 2, \, x_2(0) = 0.}
\end{cases}
\end{equation}

\subsection{Fixed Time Domain}\label{sec:numerics-fixed}
First we study the methods on a fixed time domain, $t \in [0,0.2]$, with varying time step.  We show in Figure \ref{fig:ab-am-bdf} the results from the Adams family and BDF methods.
\begin{figure}\centering
\subfloat[AB]{\includegraphics[width=0.5\columnwidth]{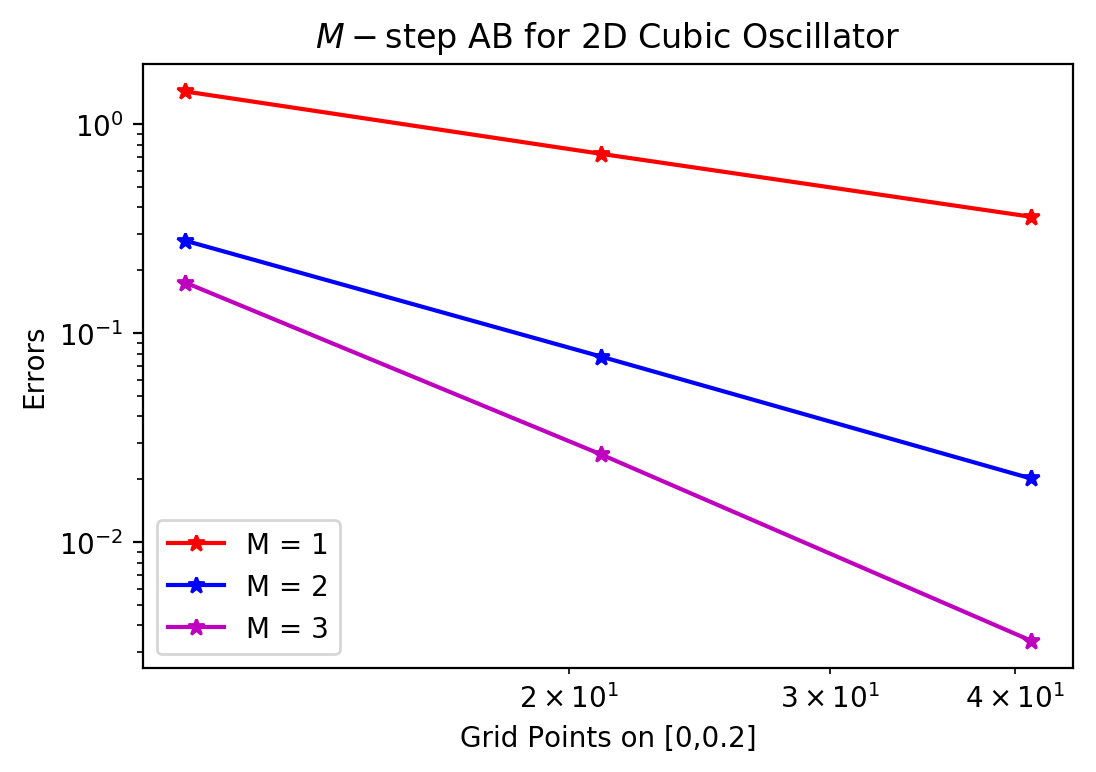}}
\subfloat[AM]{\includegraphics[width=0.5\columnwidth]{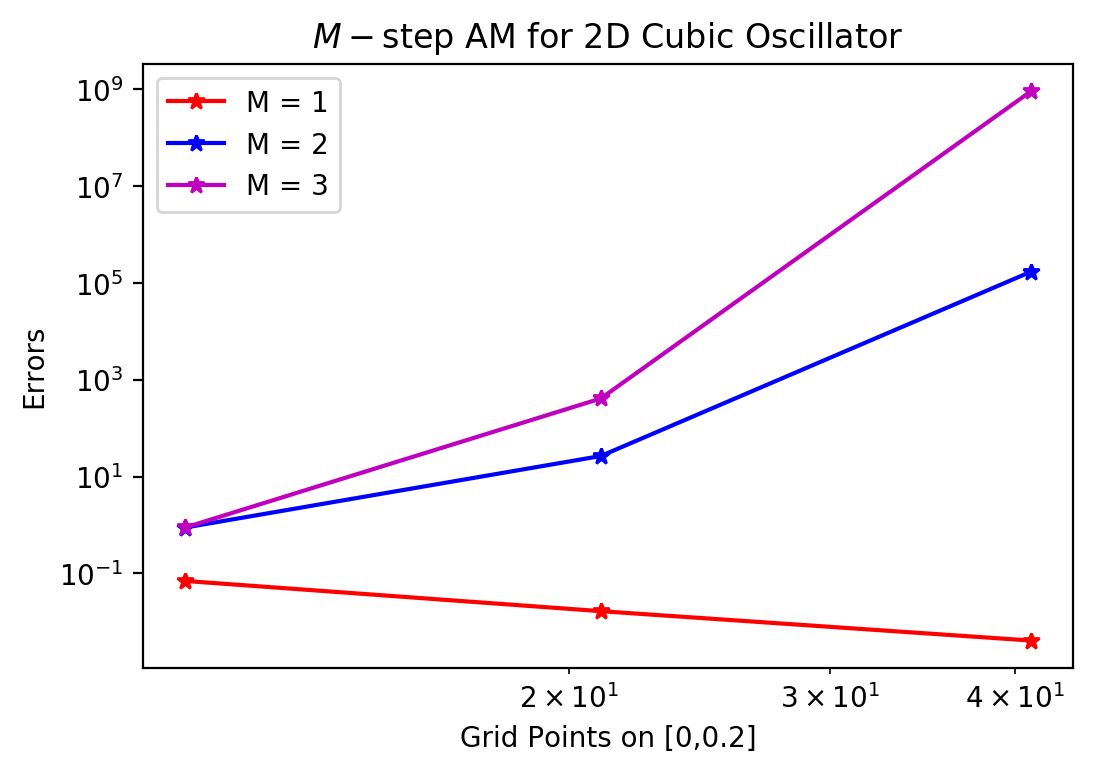}}
\\
\subfloat[BDF]{\includegraphics[width=0.5\columnwidth]{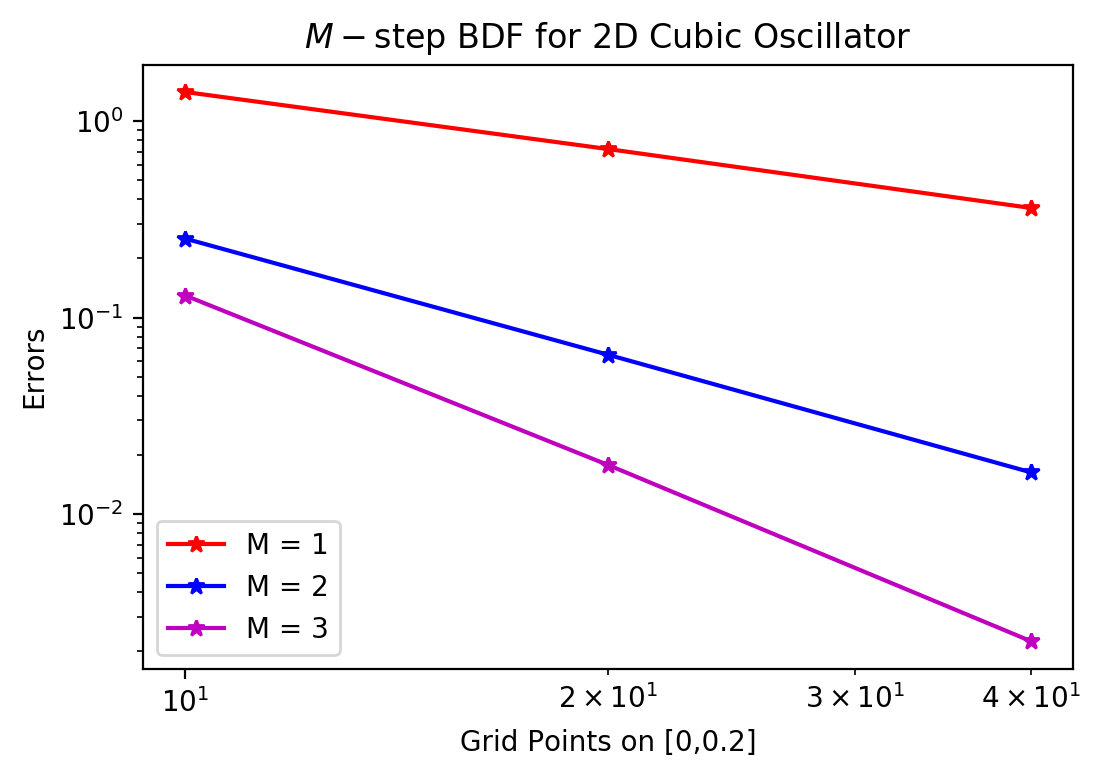}}
\subfloat[Captured Dynamics with AB-3]{\includegraphics[width=0.5\columnwidth]{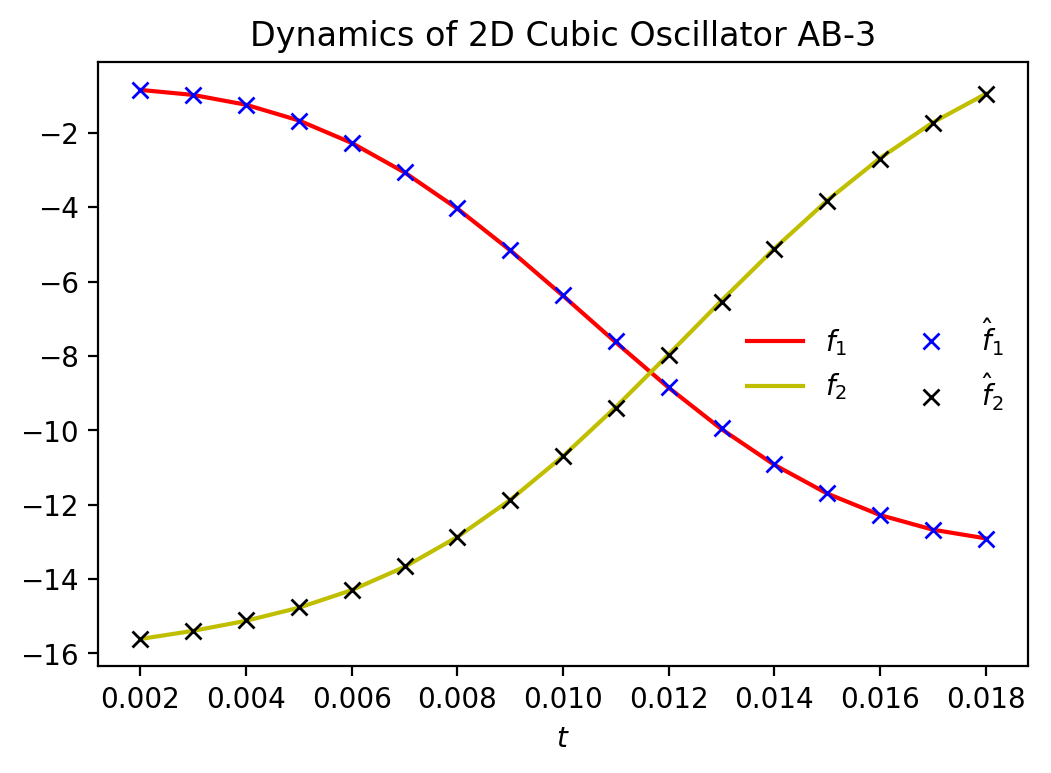}\label{ab3-dyn}}
\caption{Numerical results of the three types of schemes on the 2D cubic system \eqref{model1} on the unit time interval for different choices of $M$ and $N$.}\label{fig:ab-am-bdf}
\end{figure}
 The exact dynamics {are} computed by numerically integrating \eqref{model1} on a very refined mesh. The errors of the discovered dynamics in the $\ell^\infty-$norm are shown  in Figure \ref{fig:ab-am-bdf}  for different $M$ against different number of grid points. In addition, {Figure \ref{ab3-dyn} shows a segment of the approximated dynamics captured over the interval versus the true dynamics using a stable and 
 convergent method (AB-3) when $h=0.01$}. {In this figure, the dotted and dashed lines represent the true dynamics in the first and second coordinates, i.e. $\fb_1$ and $\fb_2$, respectively. The crosses and asterisks denote the learned dynamics in the first and second coordinates, i.e. $\hat{\fb}_1$ and $\hat{\fb}_2$, respectively. The method is able to capture the twist and intersection of the two coordinates.}
Clearly, the numerical results support the theoretical findings of this paper.


\subsection{Long Time Behavior}\label{sec:long-time-test}
{Here, we consider the problem of discovering dynamics over a changing domain $[0,T], T\gg 1.$ with fixed mesh size $h$. 
In Figure \ref{fig:long-time}, we discover the dynamics of the 2D Cubic System over specified ranges of T ($T=12.5, 25, 37.5).$  For AM (Figure \ref{fig:am-longTime}, \ref{fig:ab-longTime}, and \ref{fig:bdf-longTime}) we use  $h=0.01$ to first generate data over $[0,50]$ and then select the slice of data matching the $T$s.  AM-$M$ clearly suffers from the 
exponential error growth when $M\geq 2$, while it has a {constant error} when $M = 1$, as predicted in Section \ref{sec:long-time}.}
Meanwhile, also consistent with the analysis of Secton \ref{sec:long-time}, AB and BDF are robust for the long-time dynamics discovery -- yielding a constant error for fixed mesh as $T$ increases and a  decreasing error for larger $M$.
\begin{figure}\centering
\subfloat[AM]{\includegraphics[width=0.33\columnwidth]{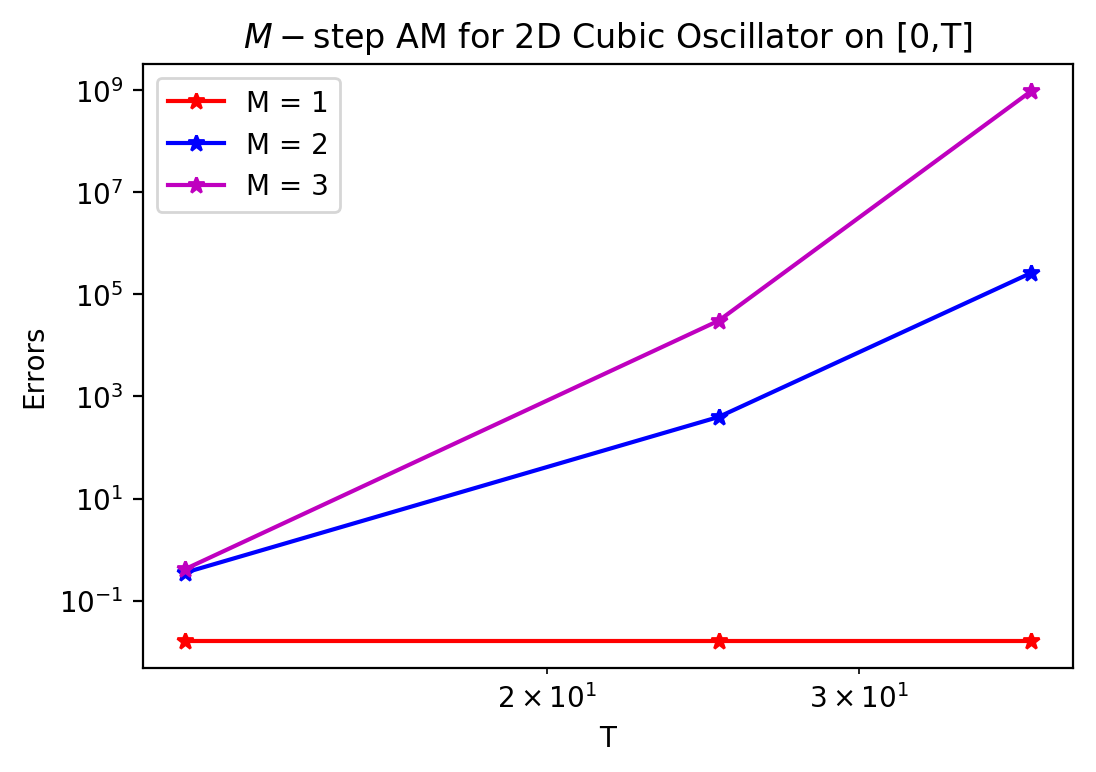}\label{fig:am-longTime}}
\subfloat[AB]{\includegraphics[width=0.33\columnwidth]{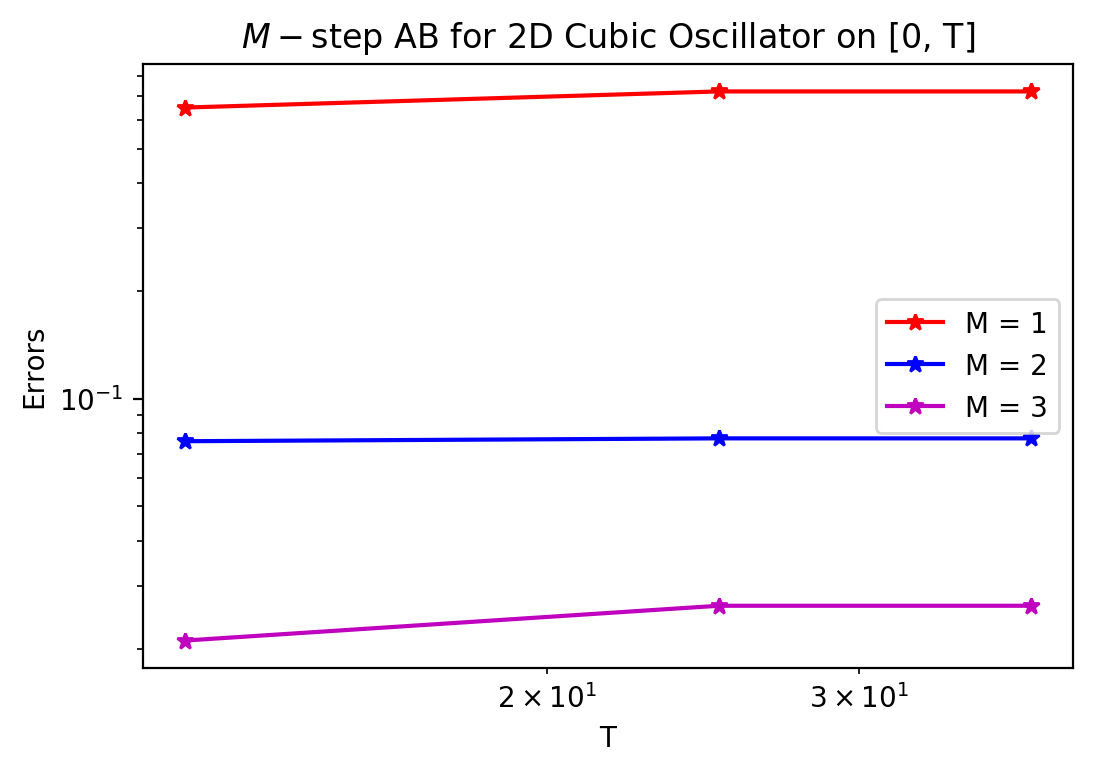}\label{fig:ab-longTime}}
\subfloat[BDF]{\includegraphics[width=0.33\columnwidth]{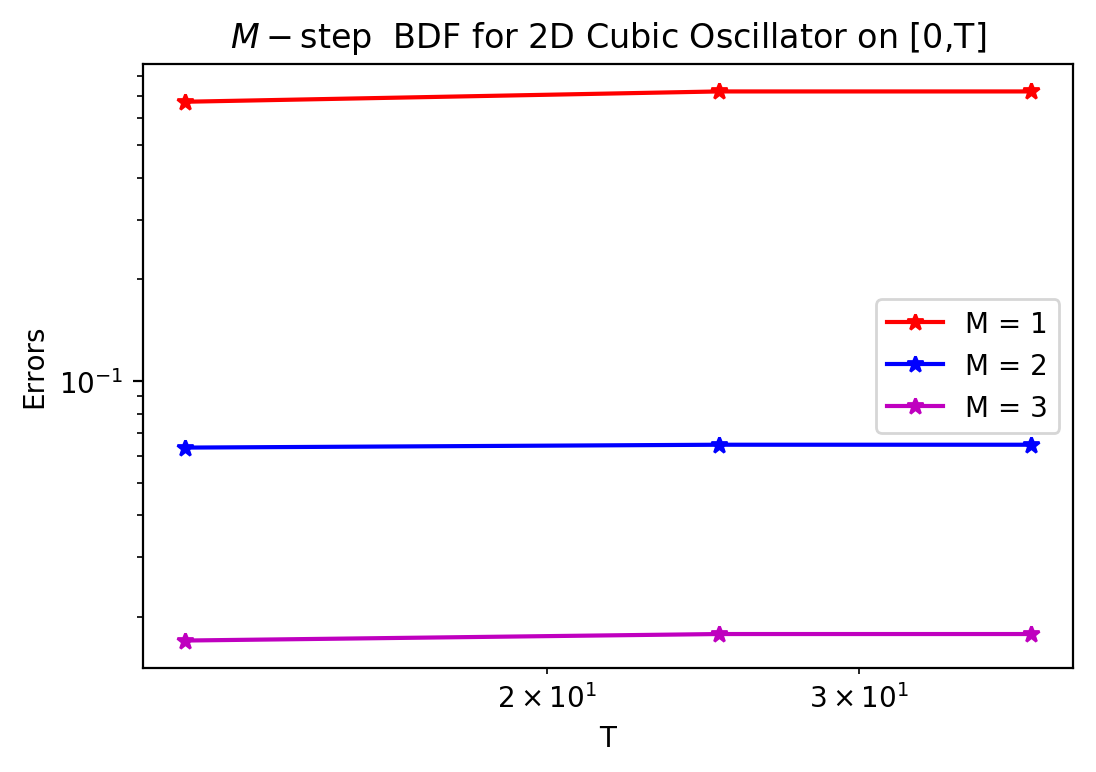}\label{fig:bdf-longTime}}
\caption{Long Time Errors for Discovery of 2D Cubic System \eqref{model1}}\label{fig:long-time}
\end{figure}


{
\begin{table}[!bp]
{\renewcommand{\arraystretch}{1.4}%
\begin{center}
\begin{tabular}{l||c|c}
\hline
Task & Integrating dynamics & Learning dynamics\\\hline 
Goal & Given {$\fb=\fb(\xb)$, find $\xb=\xb(t)$}{.} & Given {$\{\xb(t_n)\}_{n=0}^N$}, find {$\fb=\fb(\xb)$}{.}\\\hline 
Type & Forward problem & Inverse problem\\\hline 
Consistency & $\rho(1)=0, \rho'(1)=\sigma(1)$ &   $\rho(1)=0, \rho'(1)=\sigma(1)$  \\\hline 
Stability & 
{Dalhquist} root condition on $\rho$ & {New root conditions on $\sigma$  {(or $\hat{\sigma}$)} }\\\hline
Example &  {BDF-$M$ ($M\leq 6$), } AM, AB &  BDF, AM-0, AM-1, AB-$M$ ($M\leq 6$)\\\hline
\end{tabular}
\end{center}
}\caption{LMMs: similarities and differences for integrating and learning dynamics.}\label{table-LMM}
\end{table}
}

\section{Conclusions and Future Steps}\label{sec:conc}
In this paper, we extend the foundational work of solving ordinary differential equations using LMMs to the problem of dynamics discovery. We introduce refined notions of consistency, and stability, and convergence for discovery based on classical definitions, and we showed how three prominent schemes -- Adams-Bashforth, Adams-Moulton, and Backwards Differentiation Formula -- may or may not be convergent numerical methods for dynamics discovery in general.
{To do so, we first derive algebraic criteria to determine the consistency and stability of the LMM, in a spirit similar 
 to the counterpart for the classical theory. The key difference lies in the characteristic polynomial of attention; instead of the root condition 
 for the first characteristic polynomial, as classically attributed to LMMs as time integrators, stability for discovery of dynamics is attributed to root conditions on the second characteristic polynomial. While the conditions are trivial for the BDF class,
 their validity in the case of AM schemes requires the study of some new properties  of the Lagrange interpolants. The case of AB, at the present, has to be investigated computationally.
Numerical results are presented to show agreement with the theoretical findings. In conclusion, we find theoretically and numerically that the systems for  {BDF-$M$}  for all $M \in \mathbb{N}$, AB for $ 1 \leq M \leq 6$, and {AM-$M$} for $M=0$ and $1$  are and convergent, while {AB-$M$} for $7 \leq M \leq 10$ and {AM-$M$} for $M \geq 2$ are not, as summarized in Table~\ref{table-LMM}.
These conclusions are drawn provided some initial data on the dynamics  are supplied. Modifications need to be made,
as discussed in Section~\ref{sec:disc},
if other types of additional data on the dynamics are provided. }
{LMM schemes are well-studied for the forward problem in numerical analysis. 
As such tools, they can be useful to the subject of machine learning. For example, 
they can be applied to the design and training of neural networks that are seen as discrete forms of  dynamic systems \cite{weinan2017proposal,chen2018neural,sun2018stochastic}.  Different from such applications, the new study given here is motived by recent interest in using machine learning  \cite{bishop06,goodfellow2016deep,mohri18,murphy12,shalev2014} to
formalize a variety of inverse problems {such} as learning dynamics using classical discretization techniques like LMMs.  The change of the problem type from forward integration to inverse learning leads to different mathematical theory as illustrated in the Table~\ref{table-LMM}\footnote{{Note that there are several different versions of consistency and stability of LMM based dynamics discovery, which also affect the order of convergence, see discussions in 
 Theorem~\ref{thm:order}.}}.
Note that in particular, BDF provides a class of methods convergent for integrating and learning dynamics, while not all AB and AM methods can share the same conclusion. Our framework can be applied to check on other LMMs besides these examples. Furthermore, it  will be interesting to explore if there are systematic ways to generate broader classes of LMMs good for both tasks of model-based  time integration and data-driven learning.
}

{
As discussed in Section~\ref{sec:connect}, our current study assumes the best possible case that the exact states along with suitable approximations to the initial dynamics are all given, together with the assumption that the neural network representation
 can produce zero residual for the LMM dynamics.
While this setting is highly idealized,  based on the conclusions drawn, we can speculate about the impact  on the properties of stability and convergence caused by different choices of time discretization schemes for a more informed attempt at discovery of unknown dynamics in more 
practical settings. The latter leads to many interesting issues to be considered in the future.  For instance,  instead of 
assuming only data on the state with a loss function $\cT(\tilde{\xb},\tilde{\fb}, \fb_{NN})$, we may consider a 
more general loss function with data on the state and dynamics, i.e.  $\cT(\tilde{\xb}, \tilde{\fb}, \fb_{NN}, \hat{\fb},\hat{\xb})$, given by
}
$${
 \underbrace{\ell_1(\, \hat{\xb},\hat{\fb} \, )}_{\text{dynamics conformity}} \hspace{-3mm}
+\ \, \underbrace{\ell_2(\tilde{\xb},\tilde{\fb} ,\fb_{NN},\hat{\fb} ) +\ell_3(\tilde{\xb},\hat{\xb})}_{\text{data fidelity}}+\underbrace{ \mathcal{R}_1(\hat{\xb}) +
  \mathcal{R}_2(\hat{\fb})+ \mathcal{R}_3({\fb}_{NN})}_{\text{regularization}}.
  }$$
{
 For LMMs with grid functions, the loss $\ell_1$ associated with dynamics conformity comes from the discretization \eqref{dyna-eqn}, and $\tilde{\fb} \in \Gamma_h[a,b]$, the space of grid functions.}
{The total loss can be taken as an expectation over training samples
 and minimized to obtain some optimal representation of the state or dynamics.  LMNet is an example where the conformity term is minimized over parameterized neural networks  of various types, so that {$\hat{\fb} \equiv \fb_{NN}$}, as studied in \cite{qin19,mnn,perdikaris19,webster18}.}
{
Whenever the term involving the LMM residual is accounted for, the framework developed in this paper would be relevant.  For stable LMMs considered here, one may expect that it may be possible to extend the convergence results for exact and complete data if the set of neural networks can satisfy some universal approximation properties. The convergence would be expected to be in the sense of function approximations which would imply good generalization error, at least among suitable classes of smooth dynamic systems. For systems displaying chaotic behavior and sharp transitions, new ideas are likely needed in order to assure accurate discovery of the underlying complex dynamics.}

{
In this more general setting, 
neural network representations may also provide implicit regularization of the learned dynamics
so that unstable LMMs could potentially be stabilized. However,  regularization likely produces additional consistency error so the convergence has to be more carefully examined.
Moreover,  we may consider compressed representation and treat incomplete data by promoting sparsity 
or exploring the use of partial physics as regularization to achieve physics-informed and data-driven discovery of the dynamics. Finally, there are many avenues of exploration to extend the results reported here.
Some interesting topics for future studies include }
\begin{enumerate}
\item  the effects of regularization {by specifying various forms of 
 the regularization terms $\mathcal{R}_1$ and $\mathcal{R}_2$, such as those promoting smoothness, sparsity, low dimensionality, and
  extending the above tasks for study of the dynamics discovery problem with incomplete and uncertain data.}
\item {different} reduced-order models {via choice of 
 constrained representations on the dynamics or the state variables or both } {\cite{bhattacharya20,webster18}};
\item {extension of} the stability framework to incorporate other multistep {and multistage} schemes such as predictor-corrector, {Milne} and  Runge-Kutta \cite{rudy19};
\item {derivation of a general class of LMMs that are convergent for both the forward problem of time integration and the backward problem of dynamics discovery.}
\item the errors in numerically {\it integrated} states based on learned dynamics \cite{mnn};
\item distributed dynamic systems such as time-dependent PDEs and examine the additional effect due to spatial discretization;
\item {generalizing to the study of  dynamics for a suitable set of  initial conditions.}
\end{enumerate}
{Naturally, learning dynamics has strong connections to the subject of time-series prediction using deep learning \cite{connor1994recurrent,han2004prediction,karim2017lstm,tao2018hierarchical}. Our current work here may motivate further  rigorous numerical analysis studies in such a direction as well.
To conclude, we see from this study that there are many new challenges in physics-based and data-driven modeling and simulations warranting further numerical analysis research.
}

\section{Acknowledgments}\label{sec:ack}
The authors would like to thank the CM3 group at Columbia University for invigorating discussions, {Wen Ding for his stimulating suggestions, and the referees and Associate Editor of \emph{SIAM Journal of Numerical Analysis} for their valuable comments.}

\bibliographystyle{siamplain}
\bibliography{references}

\end{document}